\documentclass[letterpaper, 12 pt, onecolumn, draftcls]{ieeeconf}  
\IEEEoverridecommandlockouts                              
\usepackage{amsmath,amssymb,bbm, amsthm}
\newcommand{\C}{\mathcal{C}}
\newcommand{\G}{\mathcal{G}}
\newcommand{\V}{\mathcal{V}}
\newcommand{\E}{\mathcal{E}}
\newcommand{\I}{\mathcal{I}}
\newcommand{\K}{\mathcal{K}}

\newcommand{\R}{\mathcal{R}}
\renewcommand{\S}{\mathcal{S}}

\newcommand{\fin}{f^{\mathrm{in}}}
\newcommand{\fout}{f^{\mathrm{out}}}
\newcommand{\hatfin}{\hat{f}^{\mathrm{in}}}
\newcommand{\hatfout}{\hat{f}^{\mathrm{out}}}
\newcommand{\rhojam}{\rho^{\mathrm{jam}}}
\newcommand{\ones}{\mathbf{1}}
\newcommand{\floor}[1]{\lfloor#1\rfloor}
\newcommand{\trace}[1]{\mathrm{trace}\left\{#1\right\}}
\newcommand{\rank}[1]{\mathrm{rank}\left\{#1\right\}}
\newcommand{\virg}[1]{``#1''}
\newtheorem{remark}{Remark}
\newtheorem{lemma}{Lemma}
\newtheorem{assumption}{Assumption}
\usepackage{graphicx}
\usepackage{color}
\graphicspath{{Images/}}
\newtheoremstyle{problemstyle}  
        {3pt}                                               
        {3pt}                                               
        {\normalfont}                               
        {}                                                  
        {\bfseries\itshape}                 
        {\normalfont\bfseries:}         
        {.5em}                                          
        {}                                                  
\theoremstyle{problemstyle}
\newtheorem{problem}{Problem}
\title{\LARGE Flow and Density Reconstruction\\ and Optimal Sensor Placement\\for Road Transportation Networks}

\author{Enrico Lovisari, Carlos Canudas de Wit, and Alain Y. Kibangou
\thanks{The authors are with Univ. Grenoble Alpes, Gipsa-Lab, with CNRS, Gipsa-Lab, and with Inria Grenoble Rh\^{o}ne-Alpes, F-38000 Grenoble, France {\tt\small {enrico.lovisari, carlos.canudas-de-wit, alain.kibangou}@gipsa-lab.fr}}%
}

\begin{document}

\maketitle
\begin{abstract}

This paper addresses the two problems of flow and density reconstruction in Road Transportation Networks with heterogeneous information sources and cost effective sensor placement. Following standard macroscopic modeling approaches, the network is partitioned in cells, whose density of vehicles changes dynamically in time according to first order conservation laws. The first problem is to estimate the flow and the density of vehicles using two sources of information, namely standard fixed sensors, precise but expensive, and Floating Car Data, less precise due to low penetration rates, but already available on most of the main roads. A data fusion algorithm is proposed to merge the two sources of information for observing density and flow of vehicles. The second problem is to place the sensors in the network by trading off between cost and performance. A relaxation of the problem is proposed based on the concept of Virtual Variances. The efficiency of the proposed strategies is shown on a synthetic regular grid and in the real world scenario of Rocade Sud in Grenoble, France, a ring road 10.5 km long.
\end{abstract}

\section{Introduction}
\label{sec:introduction}

The last decades have witnessed a considerable increase of the number of vehicles especially as consequence of urbanisation in big metropolis, not matched by a comparable extension of road infrastructures. As a consequence, crucial freeways, highways and arterial roads have been steered to a state of near saturation, and experience on daily basis periods of congested traffic \cite{PapageorgiouHORMS:07}. In turn, congestion causes increased travel times and stop-and-go phenomena, leading to decreased safety, economical losses, and environmental and psychological hazards in terms of pollution and road rage \cite{BilbaoTRA:08}. Increasing road capacity by extending road infrastructures, such as construction of new arterial roads, has been the standard way to cope with congestion problems, but that is often infeasible when existing roads lie on built-in areas. Intelligent Transportation Systems (ITSs), on the contrary, exploit recent technological and theoretical advancements in distributed computation and communication, and are expected to provide robust techniques for real-time monitoring, prediction and actuation of traffic networks, and to better integrate with road and rail public transportation.

The first goal of the present paper is addressing the problem of estimating road usage in terms of density and flow of vehicles in a traffic network. The latter are commonly considered a good representation of the state of the system, providing more information than average speed alone. In particular, they are of crucial importance for 1) forecasting travel time and traffic evolution, along with historical data; 2) informing in real-time drivers about the state of the network through navigation systems; 3) providing public authorities with statistical data to monitor the state of the network and predict dangerous scenarios; 4) computing and actuating control actions through traffic lights, ramp metering and speed limits, or, in the future, lane change and semi-autonomous routing and navigation \cite{PapageorgiouIEEE:03,PapageorgiouTRR:91, PisarskiCDC:12, ComoCDC13}.

Standard devices to obtain information on the state of the network are fixed sensors such as induction loops and magnetometers. Placed over a section of road, they provide rich information on the vehicles that cross such a section over a pre-fixed period of time: 1) their number, or flow, 2) their average speed, and 3) their average density, or more precisely their occupancy (see Section~\ref{sec:model}). Current technology allows for very precise measurements, with relative errors of measured quantities against ground truth often being below $1\sim 2 \%$. However, deployment and maintenance of a sensing network requires considerable investments and human force, and consequently sensing networks are usually designed to be as sparse as possible. The second problem addressed in this paper is indeed the Optimal Sensor Placement problem, that is, positioning sensors on the cells of the network given partial information on the system and in such a way to trade off between performance and cost.
\medskip

Recently, the spread of wireless devices allows sensing and communication capabilities unforeseeable up to few years ago. Limiting the attention to traffic applications, vehicles equipped with positioning devices (such as GPS) and able to communicate with an ITS monitoring system can act as a probes in the traffic and provide Floating Car Data (FCD), namely, information on the vehicles' positions and speeds. The collected data can be used to estimate the speed in the network, thus offering a second source of information. Due to privacy reasons, single vehicles traces are usually not directly used, but rather aggregated as average speed of vehicles in segments of road. Advanced methodologies ensure fine spatial partitions of the network, with segments as short as 250 meters \cite{INRIX}.\ Compared to fixed sensors, a service based on this technology can only make use of information coming from her customers, which are a fraction of the total vehicles on the road (the \emph{penetration rate} of the system). This implies that speed measurements are less precise and flow measurement are unavailable. On the other side, since it exploits existing communication systems it is relatively inexpensive and, more important, already covers all major traffic networks.

In the first part of this paper we propose an algorithm that aims to reconstruct the traffic density and flow by fusing fixed sensors measurements and Floating Car Data. We employ a macroscopic model, partitioning the network in cells and assigning to each cell a density of vehicles. The latter evolves dynamically according to a first order mass-conservation law.

Traffic models date back to the first half of the XXth century. The most celebrated macroscopic model is the PDE based Lighthill-Whitham and Richards (LWR) model \cite{Lighthill:55}, which, based on fluid kinematics, is able to reproduce crucial phenomena such as traffic shock waves. Discretization of the LWR-PDE is not straightforward but stable numerical schemes have been proposed, the most well known being the Cell Transmission Model (CTM) \cite{Daganzo:94, Daganzo:95}. Huge efforts have been put in the last 15 years to calibrate the CTM \cite{MunozTRB:06} and to unveil its system-theoretical properties \cite{MorbidiECC:14}. Fusion of flow, density and speed measurements has also been addressed, even though mostly considering single vehicles traces. Approaches range from signal processing techniques such as the generalized Treiber-Helbing filter \cite{vanLintCACIE:10}, nonlinear versions of the Kalman filter in the context of Lagrangian sensing \cite{WorkACC:09}, and stochastic versions of the three-detector model \cite{DengTRB:13}. Recent approaches do not rely on discretization of the LWR-PDE model and allow to cast problems of estimation and control as convex problems \cite{LiTCONES:14}.

We inherit from the CTM the assumption that the inflow in a cell is a fixed linear combination of the outflows of the preceding cells. Differently from CTM, however, inflows and outflows in all the cells are estimated on the basis of the available flow measurements only. In addition, using the concept of Fundamental Diagram and the speed measurements, we compute an instantaneous (namely, only based on the latest available measurements) pseudo-measure of the density. These quantities are then the inputs for the density observer. Finally, we propose a gradient descent method to calibrate the Fundamental Diagram.

In the second part of the paper we address the problem of Optimal Sensor Placement, namely, the problem of finding the best location where to physically place sensors. This is based on trading off between two contrasting objectives: the first, to maximize the performance of state reconstruction; the second, to minimize the total economic cost of the network.

The performance of the state reconstruction is usually related to the ability to properly estimate the density of vehicles in the roads. Unfortunately, nonlinearity and complexity of traffic systems make it hard to evaluate the performance of nonlinear observers. In order to simplify the setting, we consider the related problem of reconstruction of flows in a static setting. In particular, we consider as performance metric the error covariance of an estimator of the cumulative flows in the network over a long period of time. The Optimal Sensor Placement problem can be then seen as trading off between the performance of such a flow estimator, and a cost that depends on the dimension of the sensing network. Since this is a combinatorial problem, we relax it using a method that we call Virtual Variance algorithm, based on the idea to associate to each sensor a virtual variance which is large when the sensor is not relevant for good reconstruction of the flow vector. The only input that the algorithm needs is the matrix of splitting ratios, that prescribes how vehicles split at each junction, and the nominal variance of each sensor. Furthermore, we discuss in detail two extensions of the proposed algorithm dealing with important scenarios. In the first, Optimal Sensor Placement with geographical constraints, we address the scenario in which sensors cannot be placed in a subset of cells of the network. In the second, Optimal Sensor Placement with Number of Sensors constraints, we deal with the case in which the maximum number of sensors is pre-specified, for example due to budget limitations.

Optimal Sensor Placement is an ubiquitous problem that has received a high degree of attention in several communities due to its importance for netwFork design. In Transportation Systems, it is of interest both in the dual-problem of best placement of hubs for cost-efficient transportation of goods and people \cite{ShahabiTRPE:14} and Origin-Destination coverage \cite{EhlertTRB:06, HuTRB:14}. In these works, and differently from the present paper, the problem is cast as a mixed integer problem which corresponds to determine the minimal set of locations from which the flows on the whole network can be determined, and sensor measurements are assumed to be perfect.


To summarize, the contributions of this paper are the following: 1) we formalize the problem and we design an easily implementable approach to data fusion of fixed sensors measurements and Floating Car Data; 2) we propose a gradient descent calibration algorithm of the underlying macroscopic model; 3) we formalize the problem of Optimal Sensor Placement in terms of positions of sensors in a network when sensors are noisy and we provide an approximate solution using the concept of Virtual Variances; 4) we show the prowess of the devised Optimal Placement procedure on a regular grid, for which we offer a comparison between the solution found with our approach and the true optimal placement, found by exhaustive search; 5) we illustrate the performance of the Optimal Placement procedure and of the Reconstruction algorithm through extensive numerical experiments in the real-world scenario of Grenoble Traffic Lab (GTL) \cite{CanudasIEEECS:15}, a sensing network deployed along the freeway ``Rocade Sud'' in Grenoble, France, with FCD provided by INRIX, one of the most well known traffic solutions companies.

The remainder of the paper is organized as follows: after setting up the notation, Section~\ref{sec:model} describes the model for a Road Transportation Network. Section~\ref{sec:densityProblemFormulation} formulates the problem of flow and density reconstruction and describes the proposed nonlinear observer, while the problem of optimal sensor placement is formulated and a solution based on the heuristic Virtual Variance algorithm is presented in Section~\ref{sec:optimalProblemFormulation}. Finally, Section~\ref{sec:numericalExperiments} illustrates the solutions on a regular grid and on the real world scenario of the freeway Rocade Sud in Grenoble, and Section~\ref{sec:conclusions} draws the conclusions and provides several future research directions.

\subsection{Notation}
\label{subsec:notation}

The symbols $\mathbb{R}^n$, $\mathbb{R}_+^n$ and $\mathbb{R}^{n \times m}$ denote the sets of real valued vectors of dimension $n$, of positive real valued vectors of dimension $n$, and of real valued matrices of dimension $n\times m$, respectively. The symbol $\mathbb{R}^\mathcal{A}$ ($\mathbb{R}_+^\mathcal{A}$) with $\mathcal{A}$ a finite set is to be interpreted as the set of real vectors (real positive vectors) indexed by elements of $\mathcal{A}$. The transpose of $A\in\mathbb{R}^{n\times m}$ is denoted $A^T$. For a vector $x\in\mathbb{R}^n$, $x\geq 0$ is meant component-wise. The symbol $I$ denotes the identity matrix of suitable dimensions. $|\mathcal{A}|$ denotes the cardinality of the set $\mathcal{A}$.

A graph $\G$ is a pair $(\V, \E)$ where $\V$ is called the set of nodes and $\E$ the set of edges. The functions $t:\,\E \to \V$ and $h:\,\E\to\V$ tell, for each edge $e$, which node is its tail and which node is its head, namely $e = (h(e), t(e))$. For $e\in\E$, denote by $\E_e^+:= \{j\in\E:\,h(e) = t(j)\}$ and $\E_e^-:= \{j\in\E:\,h(j) = t(e)\}$ the set of edges that follow or precede $e$, respectively. A path of length $n\geq2$ is a sequence of edges $e_1, \dots, e_n$ such that $e_{i+1} \in \E_{e_i}^+$ for all $i = 1, \dots, n-1$. A path of length $1$ is a path made of a single link. The matrix $L \in \mathbb{R}^{\V\times\V}$ is a sublaplacian of $\G$ if $L_{ej} > 0$ only if $(e,j) \in \E$, $e\neq j$, and $\sum_jL_{ej}\leq 0$.

\section{Road Transportation System Model}
\label{sec:model}

We adopt a macroscopic approach by partitioning the lanes of the roads in a traffic network in cells. Cells that lie on the same section of a road and on different lanes are said to be parallel one each other. We interpret each cell as an edge $e\in\E$ in a graph $\G = (\V, \E)$, which models the whole network. Here, nodes $v\in\V$ represent junctions, at which many roads intersect, or sections of a road. Among the set of cells $\E$, we denote by $\R^i$ and $\R^o$ the set of onramps and offramps, respectively. In this paper, we shall call onramp (offramp) any gate, be it a real ramp, a connector, a secondary road, etc, which lets vehicles enter into (exit from) the network.

We make the following connectivity assumption, which formalizes the mild requirement that any cell can be reached from an onramp and that vehicles from any cell can exit from the network.

\begin{assumption}
For any cell $e$ in $\E$, there is at least one onramp $j \in \R^i$ and one offramp $k\in\R^o$ such that $e$ is an edge of the path from $j$ to $k$. 
\end{assumption}

Time is discrete and slotted in intervals of duration $T > 0$. On each cell $e\in\E$, denote by $\rho_e(t)$ the density of vehicles, in number of vehicles per km\footnote{While we employ SI units for simplicity and for coherence with the data in the Grenoble Traffic Lab, the presentation of our results would obviously be unchanged if other systems of measurements, such as the imperial system, were used.}, at time $t$, and let $\rho(t) = [\rho_e(t)]_{e\in\E}$, which we call the \emph{state} of the network. The density of vehicles in a cell changes dynamically in time according to the following mass-conservation first-order model
	\begin{equation}
		\label{eq:model}
		\rho_e(t+1) = \rho_e(t) + \frac{1}{\ell_e}(\fin_e(t) - \fout_e(t)), \qquad \forall e\in\E
	\end{equation}
where $\ell_e$ is the length of cell $e$, and $\fin_e(t)$ and $\fout_e(t)$ are the inflow and the outflow at cell $e$ during the $t$-th time slot. 

To relate inflows and outflows we resort to the standard concept of splitting ratios. Indeed, denote by $R_{ek}\geq0$ the fraction of vehicles that turn into cell $k$ when they exit from cell $e$, which is the splitting ratio of $e$ towards $k$. Clearly, $R_{ek} = 0$ if $e$ and $k$ are not consecutive, and $\sum_kR_{ek} = 1$, if $e \not\in\R^o$. From this moment on, we make the following assumption.

\begin{assumption}
The set of splitting ratios $\{R_{ej}\}_{(e,j)\in\E\times\E}$ is known.
\end{assumption}

Splitting ratios establish the relation $\fin_e(t) = \sum_{j\in\E}R_{je}\fout_j(t)$, for all $t\geq0$ and for any cell $e\not\in\R^i$, while $\fin_e(t) = \lambda_e(t)$ for $e\in\R^i$, where $\lambda_e(t)$ is an exogenous external inflow. We set for sake of convenience $\lambda_e(t) = 0$ for $e\not\in\R^i$. By stacking inflows and outflows into vectors $\fin(t)$ and $\fout(t)$, respectively, we can rewrite the previous relation in matrix form as
	\begin{equation}
		\label{eq:fin}
		\fin(t) = R^T\fout(t)
	\end{equation}
where the matrix $R = [R_{ej}]_{e\in\E, j\in\E\setminus\R^i}$ is the \emph{matrix of splitting ratios}. In the present paper it is assumed that the matrix of splitting ratios is fixed, predetermined, and known. Its calibration, closely related to the estimation of Origin-Destination pairs, can be performed on single-lane freeways with onramps and offramps by taking ratios of flows on main line and ramps \cite{MunozTRB:06}. We plan to extend this setting to networks by casting the problem as an optimization problem in future research.

\subsection{On modelling of cell flows}
\label{subsec:modelling_flow}

Macroscopic models such as the CTM postulate that the flow $\fout_e(t)$ that exists from cell $e$ at time $t$ is a deterministic function of the densities in the cells around $e$. In the simplest case, where $e$ and $j$ are two consecutive cells and $\E_e^+ = \{j\}$ and $\E_j^- = \{e\}$, then in the CTM
	$$
		\fout_e(t) = \min\{d_e(\rho_e(t)), s_j(\rho_j(t))\}
	$$	
where $d_e(\rho_e(t))$ and $s_j(\rho_j(t))$ are the \emph{demand} of cell $e$ and the \emph{supply} of cell $j$, and represent the maximum outflow from $e$ and the maximum inflow into cell $j$, respectively. The resulting system is a Godunov scheme for discretization of the LWR-PDE model, and can be extended to the network case in various ways \cite{Daganzo:95, LovisariDensityCDC:15, CooganACC14}. While these models reproduce important phenomena that must be taken into account when modelling and controlling traffic networks, such as the movement of shockwaves, each of them can only partially represent traffic dynamics in networks. 

For this reason, we will avoid to explicitly model the relation between flows and density, and we will limit to the standard 
	\begin{equation}
		\label{eq:fout}
		\fout_e = \rho_ev_e, \forall e\in\E
	\end{equation}
namely that the volume that exists from a cell in a period is proportional to the density of the vehicles in the cell, and to their speed. Again, we leave unmodelled the relation between these two quantities, because, as we will make clear in the following, we assume to have a direct measurement of the average speed in each cell (in the form of Floating Car Data). In conclusion, we consider from now on the dynamics of the real system to be dictated by Eqs.~\eqref{eq:model}-\eqref{eq:fout}, where $v_e$, $e\in\E$, is an unmodelled quantity which depends on the local state of the network.
\medskip

While in the model of our network we do not use an explicit relation between flows and densities, we shall use it for data fusion and estimation purposes. In particular, we write 
\begin{equation}
		\label{eq:foutEstimation}
		\varphi_e = \varphi_e(\rho_e), \forall e\in\E
	\end{equation}
where $\varphi_e$ is the flow of vehicles \emph{at the sensor locations}, which we shall always assume to the at the end of the cell. The graph of the function $\varphi_e(\cdot)$, which is the \emph{Fundamental Diagram} on cell $e$, is a concave function with $\varphi(0) = \varphi(\rhojam) = 0$, where $\rhojam$ is the jam density, at which vehicles are too close one each other to move (values for $\rhojam$ vary from $150$ to $300$ vehicles per km).

\subsection{Available measurements}
\label{subsec:available_meas}

In this paper we consider two heterogeneous sources of information: flow and density measurements from sensors and Floating Car Data. 

\subsubsection{Flow and density measurements} 

Standard measurement devices for traffic are loop detectors or magnetometers, radar traffic detectors, or video detection systems. They are positioned at fixed and predefined positions in the network, monitor a section of road, and are able to detect and assign a timestamp to the event \virg{a vehicle crossed the section}, and information is then aggregated in time slots. For sake of simplicity and without loss of generality, we assume that such time slots correspond to the time discretization of the system \eqref{eq:model}. As such, measurements of the following two quantities are then available at every time slot of duration $T$: 
	\begin{itemize}
		\item Flow of vehicles $\varphi_e(t)$, by counting the number of vehicles crossing the section during the $t$-th time slot;
		\item Density of vehicles $\rho_e(t)$ over the section. The quantity that is actually measured by the aforementioned devices is the occupancy of vehicles $o_e(t)$, defined as the percentage of time any vehicle was standing over the section in the $t$-th time slot. However, occupancy is approximatively related to the density of vehicles by the relation $\rho_e \approx \frac{o_e}{100\ell_{\mathrm{ave}}}$, where $\ell_{\mathrm{ave}}$ is the average length of a vehicle in km. For this reason, we shall assume from now on that devices can measure density directly.
	\end{itemize}

Measurements hardly correspond to the real value of the quantities that they represent, with sources of noise ranging from temporary inability to detect changes of the magnetic field, too fast or too slow vehicles, blurred videos, etc. We adopt a simple additive noise model
	\begin{equation}
		\label{eq:flowDensityMeas}
		\begin{array}{ll}
			\varphi^m_e(t) 	 = \varphi_e(t) + \omega^{\varphi}_e(t), 	&e\in \E^m	\\
			\rho^m_e(t) 	 = \rho_e(t) +  \omega^{\rho}_e(t), 		& e\in \E^m	
		\end{array}
	\end{equation}
where $\varphi^m_e(t)$ and $\rho^m_e(t)$ are flow and density estimates at time $t$, and $\omega^{\varphi}_e(t)$ and $\omega^{\rho}_e(t)$ are measurement errors whose stochastic properties depend on the performance of the sensor as well as on road and weather conditions, and $\E^m \subseteq \E$ is the set of cells equipped with sensors. Due to installation and maintenance costs, usually $|\E^m| << |\E|$. 

\begin{remark}
It should be mentioned that measurement systems based on magnetometers can be used to measure the average speed of vehicles crossing the section they monitor, in addition to flow and density, simply by deploying them in pairs monitoring sections that are at a fixed and known distance. The two consecutive instants at which the same vehicle crosses the two sections provides then a measurement of its speed. This type of installation is however not standard and more expensive due to additional hardware and software. For this reason, and to show that our approach does not need this additional information, we shall assume that no speed measurement from static sensors is available.
\end{remark}

\subsubsection{Speed measurements}

As already mentioned, recent technological advancements provide ITSs and companies selling traffic solutions with speed measurements in the form of Floating Car Data (FCD), averaged over predefined sections of roads for privacy reasons. While some classes of vehicles, such as taxi and buses, can indeed be traced, we decided not to use such an additional information in this paper, and leave the possibility to use it as future research.

Floating Car Data are less expensive and require less maintenance effort with respect to fixed sensors since they exploit existing communication architectures. For the same reason, they are already accessible almost everywhere once a data collecting mechanism is deployed. We partition thus the network in segments, let $\S$ be the set of all segments, and we assume that a measurement of average speed is available for every segment $s\in\S$.

\begin{figure}
\centering 
\begin{tabular}{c}
\includegraphics[width=0.46\textwidth]{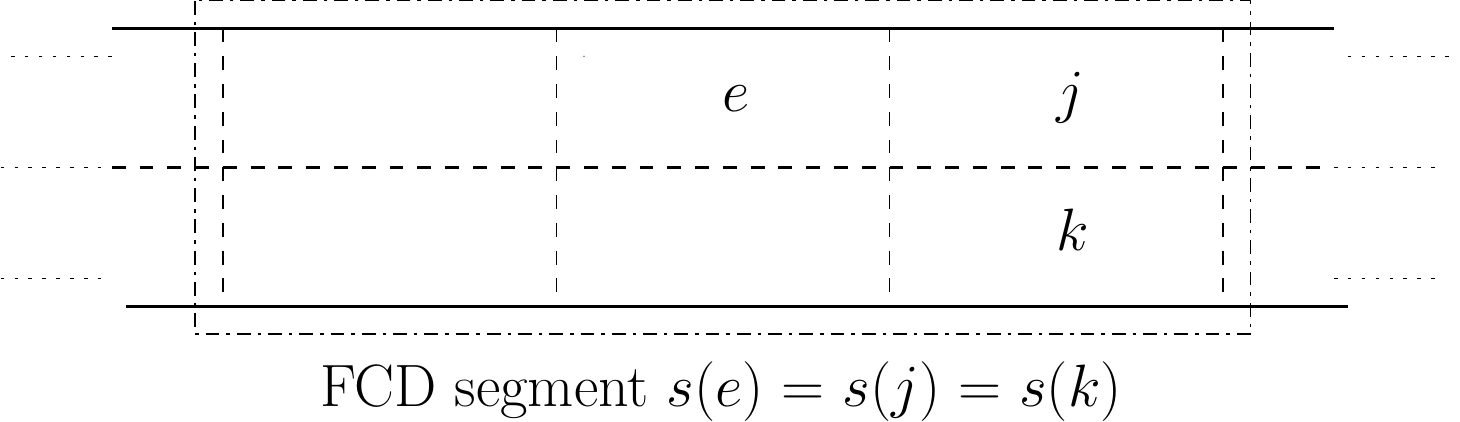}
\end{tabular}
\caption{A stretch of road partitioned in cells and FCD segments. Splitting ratios are shown from a cell $e$ to following cells, $j$ and $k$. A FCD segment including, among others, cells $e$, $j$ and $k$, is also shown.}
\label{figure:ExampleCellSegment}
\end{figure}

\begin{figure}
\centering
\begin{tabular}{cc}
\includegraphics[height=6cm]{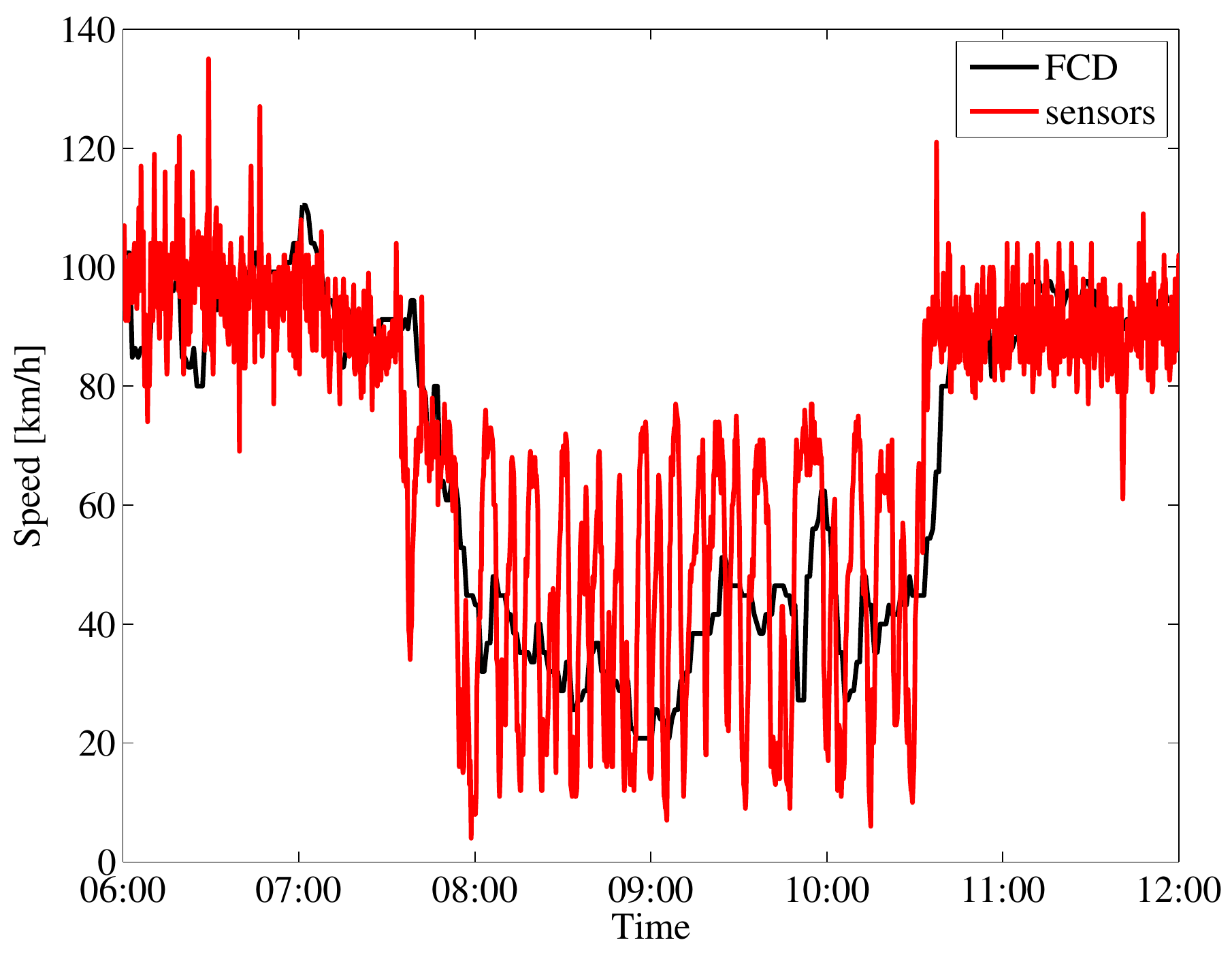}&
\includegraphics[height=6cm]{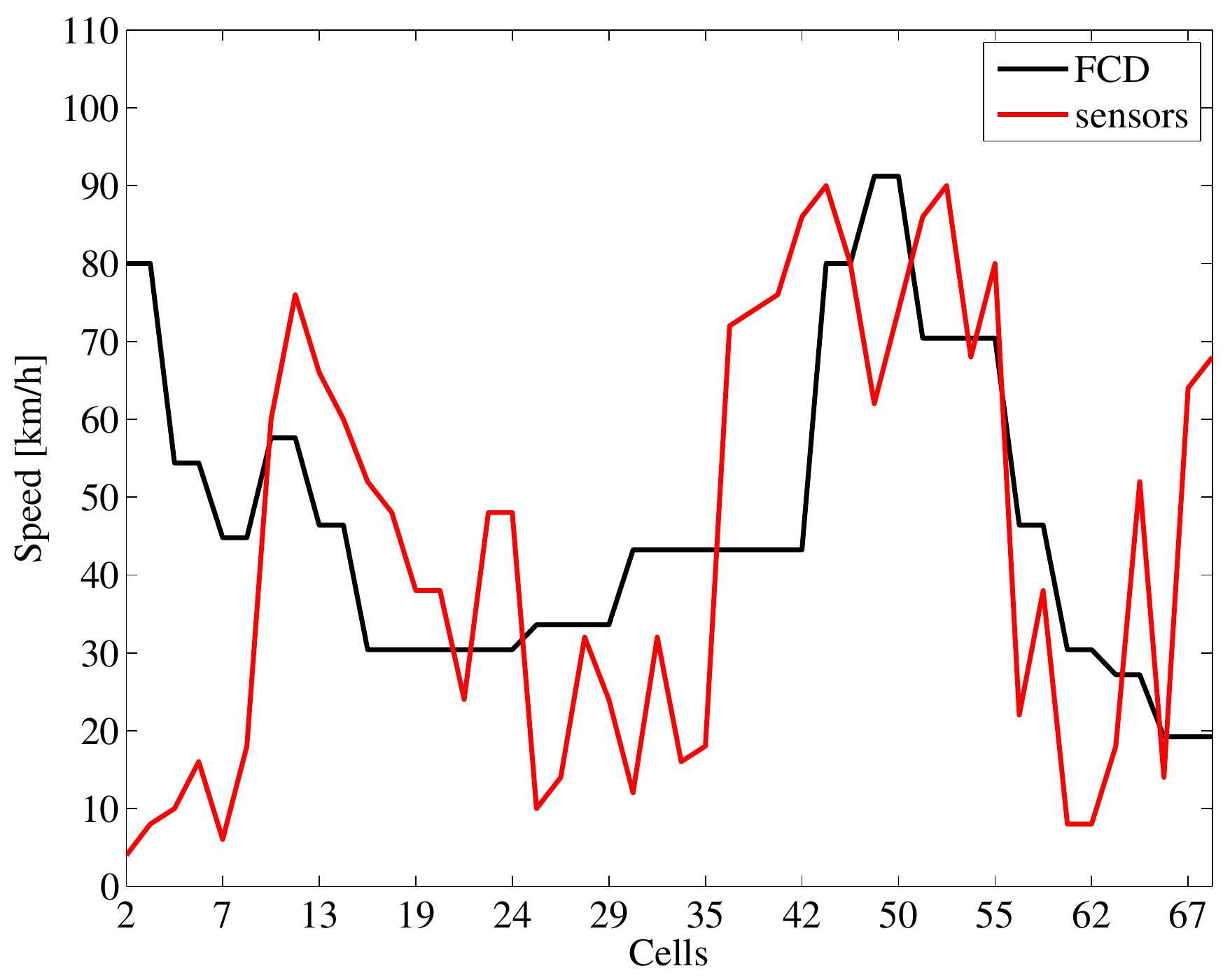}
\end{tabular}
\caption{Comparison of speed measurements from sensors and Floating Car Data.}
\label{fig:example_GTL_vs_FCD}
\end{figure}

Despite such advantages, FCD also have drawbacks. Aside from the already mentioned possibly low penetration rate, information provided via FCD is usually averaged over a relatively long period of time. As an example, within the Grenoble Traffic Lab fixed sensors yield flow and density (and speed) measurements every $T = 15$ seconds. The FCD provided by INRIX are instead aggregated per minute, with standard practice ranging between $5$ and $10$ minutes. A comparison between sensors speed measurements from the GTL and FCD is provided in Figure~\ref{fig:example_GTL_vs_FCD}. On the left panel, a comparison of FCD measurement and average of GTL for slow and fast lane at location Taillat on the Rocade Sud in Grenoble from 06:00:00 to 12:00:00 on April 24th, 2014 (see Paragraph~\ref{subsubsec:RocadeSud} for details), which clearly illustrates the high measurement rate of sensors and the averaging effect of FCD, resulting in a smoother signal. On the right panel, a comparison of the two for all the cells on the main line at 08:00:00, showing that on average the measurements are in agreement.

To formalize the scenario, we consider new speed aggregate data to be available every $N$ time instants, $N = kT$ for some $k\geq1$, namely, at times $N$, $2N$, $3N$, $\dots$, corresponding to the average speed in the periods $[0, N-1], [N, 2N-1]$, etc., respectively. As such, speed measurements can be formally written as
	\begin{equation}
		\label{eq:speedMeas}
		v^{\mathrm{FCD}}_e(t) = 
			\begin{cases}
				v^{\mathrm{ff}}_e, &t \in [0, N-1]	\\
				v^{\mathrm{FCD}}_{s(e)}(k),
				&t \in [kN, (k+1)N-1]
			\end{cases}
	\end{equation}
where
	\begin{itemize}
		\item $v^{\mathrm{ff}}_e > 0$ is the freeflow speed on cell $e$, namely, the speed of vehicles in low density regime;
		\item $v^{\mathrm{FCD}}_{s(e)}(k)$ is given by 
				$$
					v^{\mathrm{FCD}}_{s(e)}(k) = \frac{1}{N|s(e)|}\sum_{j\in s(e), \tau\in\I_t}v_j(\tau) + \omega^{\mathrm{FCD}}_{s(e)}(k)
				$$	
			where $s(e)$ denotes the segment of which $e$ is one of the cells (see Figure~\ref{figure:ExampleCellSegment}), $\omega^{\mathrm{FCD}}_{s(e)}(k)$ is a measurement error whose stochastic properties depend on the performance of the sensor as well as on road and weather conditions, and $\I_t = \{\tau:\,\floor{\frac{t}{N}}-1 \leq \frac{\tau}{N} < \floor{\frac{t}{N}}\}$.
\end{itemize}

\section{Density Reconstruction}
\label{sec:densityProblemFormulation}

The first problem we address in this paper is Density Reconstruction on the basis of heterogeneous sources of information. In particular, our aim is to build an observer for the densities of vehicles in all cells of the network given static sensor measurements and Floating Car Data.

We start by observing that Eqs.~\eqref{eq:model}-\eqref{eq:fin} cannot be directly used to observe the system except for the ideal scenario in which ideal measurements of the outflows $\fout_e(t)$, for all $e\in\E$ and for all times $t\geq0$, and of the initial conditions of the system, are available. Such a naive observer would however be very sensible to noise, as notice that while errors in the initial conditions correspond to offsets during the evolution of the system, noises in the flow measurements are integrated by the system's dynamics, thus possibly producing unbounded and/or unrealistic results. Since real systems are never error free, Eqs.~\eqref{eq:model}-\eqref{eq:fin} cannot be directly used to observe the system. 

We solve this difficulty by considering the following standard Luenberger-like observer
	\begin{equation}
		\label{eq:observer}
			\begin{cases}
				\hat{\rho}_e(t+1) = \hat{\rho}_e(t) + \frac{1}{\ell_e}(\hatfin_e(t) - \hat{\fout}_e(t)) + \kappa (\tilde{\rho}_e(t) - \hat{\rho}_e(t))	\\
				\hatfin_e(t) = \hatfin_e(\varphi^m(t))	\\
				\hatfout_e(t) = \hatfout_e(\varphi^m(t)))	\\
				\tilde{\rho}_e(t) = \tilde{\rho}_e(\varphi^m(t), v^m(t))
			\end{cases} \forall e\in\E
	\end{equation}
where
	\begin{itemize}
		\item $\hat{\rho}_e(t)$ is the estimate of the density on cell $e$ at time $t$;
		\item $\hatfin_e(t)$, $\hatfout_e(t)$ are estimates, based on the flow measurements, of inflow and outflow in cell $e$ at time $t$;
		\item $\tilde{\rho}_e(t)$ is an density pseudo-measure, based on flow and speed measurements, of the density on cell $e$ at time $t$;
		\item $\kappa$ is a tunable gain trading off between flow and density pseudo-measure;
		\item $\varphi^m(t) = [\varphi^m_e(t)]_{e\in\E_m}$ and $v^{\mathrm{FCD}}(t) = [v^{\mathrm{FCD}}_e(t)]_{e\in\E}$ are the stacked versions of flow and speed measurements.
	\end{itemize}

We are interested in the following:

\begin{problem}[Flow and Density Reconstruction using Heterogeneous Sources]
\label{pr:density}

Design the maps
	\begin{align*}
		\hatfin &= \{\hatfin_e\}_{e\in\E}:\,\mathbb{R}^{\E^m}_+ \to \mathbb{R}^\E\\
		\hatfout &= \{\hatfout_e\}_{e\in\E}:\,\mathbb{R}^{\E^m}_+ \to \mathbb{R}^\E\\
		\tilde{\rho}&=\{\tilde{\rho}_e\}_{e\in\E}:\, \mathbb{R}^{\E^m}_+\times\mathbb{R}^{\E}_+ \to \mathbb{R}^\E
	\end{align*}
to minimize the absolute errors with respect to real flows and densities
	\begin{equation}
		\label{eq:metrics_reconstruction}
		\begin{split}		
		a^{\rho}(t, e) = \left|\hat{\rho}_e(t) - \rho_e(t)\right|\\
		a^{\varphi}(t, e) = \left|\hat{\fout}_e(t) - \fout_e(t)\right|	
		\end{split}
	\end{equation}
\end{problem}

\subsection{A nonlinear observer for traffic networks}
\label{subsec:observer}

In this subsection we describe the proposed solution to Problem~\ref{pr:density}. It consists in an offline calibration procedure and an online filtering step.

\subsubsection{Offline calibration}
\label{subsubsec:offline_calibration}

This paragraph is devoted to providing a solution for calibrating the Fundamental Diagram. We employ a gradient descent strategy, similarly to \cite{QuTRB15}, for which we do not require a CTM formulation as in \cite{MunozACC04}.

Recall that the Fundamental Diagram is the graph of the function $\varphi_e(\cdot)$ that is the flow of vehicles at the point where sensors are placed on cell $e$ and whose argument is the density of vehicles on cell $e$. As such, it can be only estimated on cells $e\in\E^m$, where measurements of flow and density are available. For a cell $j \in \E\setminus\E^m$, we assume that the Fundamental Diagram can be estimated by extending by linear interpolation the parameters on the cells in $\E^m$ that are close to $j$, with coefficients that depend on the mutual distance between the cells and on the type of road $j$ belongs to.

The profile of Fundamental Diagram that we consider in the present paper is the following
	$$
		\varphi_e(\rho)=
			\begin{cases}
				v^{\mathrm{ff}}_e \rho	,&\rho \leq \rho^c_e	\\
				a_e \rho^2 + b_e\rho + c_e, &\rho > \rho^c_e
			\end{cases}
	$$
where 
	\begin{itemize}
		\item $\rho_e^c$ is the \emph{critical density}. It partitions the set of densities $[0, \rho^{\mathrm{jam}}]$ into the \emph{freeflow} low-density region $[0, \rho_e^c)$, in which the mutual influence of vehicles is small, from the high-density \emph{congested} region $(\rho_e^c, \rho_e^{\mathrm{jam}}]$, in which speed decreases with density due to interaction of close vehicles;
		\item $\rho_e^{\mathrm{jam}}$ is the jam density, at which vehicles are so close one each other that they are unable to travel;
		\item $v^{\mathrm{ff}}_e > 0$ is the freeflow speed on cell $e$; the value $C_e = v^{\mathrm{ff}}_e\rho_e^c$ is the \emph{capacity} of the section or road, namely, the maximum number of vehicles that can flow through it during a period $T$;
		\item we assume that the Fundamental Diagram congested region is a convex quadratic function of the density. The following relations among the parameters $a_e$, $b_e$ and $c_e$ hold for consistency
			$$
				\begin{cases}
					a_e\left(\rho_e^c\right)^2 + b_e\rho_e^c + c_e = v^{\mathrm{ff}}_e \rho_e^c	 \\
					a_e\left(\rho_e^{\mathrm{jam}}\right)^2 + b_e\rho_e^{\mathrm{jam}} + c_e = 0 \\
					a_e \geq 0
				\end{cases}
			$$
	\end{itemize}
	
\begin{remark}
We defined $\varphi_e(\cdot)$ to be the flow of vehicles through a section during a sample time $T$. As such, the units of the speed $v^{\mathrm{ff}}_e$ are km per $T$. 
\end{remark}	

\begin{remark}
A standard choice for the Fundamental Diagram is the a triangular Fundamental Diagram, for which in the congested region
	$$
		\varphi_e(\rho_e) = \omega_e(\rho_e^{\mathrm{jam}} - \rho_e)
	$$
where $\omega_e$ is the \emph{wave speed} at section $e$. Clearly, our model recovers the latter with $a_e = 0$, $b_e = -\omega_e$, and $c_e = \omega_e\rho_e^{\mathrm{jam}}$. The choice of a quadratic Fundamental Diagram has been driven by the empirical observation, based on data on our experimental setting, that the triangular diagram tends to overestimate the flow in congestion, as it will be shown in Section~\ref{sec:numericalExperiments}.

An alternative appealing solution which fits our data is the inverted-$\lambda$ fundamental diagram \cite{HallTRA:86}. However, the number of parameters to be estimated is higher in the latter case, and the resulting model is more complex as it involves hysteresis. 

Motivated by the well known issue that deterministic Fundamental Diagrams are in any case only a rough approximation of the relation between flow and density, we chose the quadratic profile because of it is simple to calibrate and to use.
\end{remark}

We describe now the procedure for the calibration of the Fundamental Diagram. Let $e\in\E^m$ and let $\{(\rho_{e,k}^m, \varphi_{e,k}^m)\}_{k\in\K}$, $\K = \{1, \dots, K\}$, the set of $K$ density and flow measurements used as learning set and obtained via the fixed sensor on cell $e$. For sake of notation, and since all variables refer to cell $e$ only, let us write from now on $\rho_{k}^m$ and $\varphi_{k}^m$ instead of $\rho_{e,k}^m$ and $\varphi_{e,k}^m$, and the same for the parameters of the Fundamental Diagram. 

The proposed calibration procedure requires two steps
	\begin{itemize}
		\item Estimation of $\rho^c$ and $C = v^{\mathrm{ff}} \rho^c$: the first step consists in estimating the critical density and the capacity $C$ of the cell. We consider the standard least square estimation, which results into the following non-linear and non-convex minimization problem: given the set of measures $\{(\rho_k^m, \varphi_k^m)\}_{k\in\K}$, $\K = \{1, \dots, K\}$, solve
			\begin{equation}
				\label{eq:FDCalibrationOriginal}
				\begin{array}{ll}
				\min_{(\rho^c, C)}  & V_{(\rho^c, C)} = \frac{1}{2}\sum_{k=1}^K(\varphi_k - \varphi_{(\rho^c,C)}(\rho_k))^2	\\
				\mathrm{s.t.}		& 0 < \rho^c < \rhojam	\\
									& C > 0	\\
									& \varphi_{(\rho^c,C)}(x) = 
										\begin{cases}
											\frac{C}{\rho^c}x, & x \leq \rho^c	\\
											\frac{C(\rhojam-x)}{\rhojam - \rho^c}, & x > \rho^c	\\											
										\end{cases}
				\end{array}
			\end{equation}
			
We aim to solve \eqref{eq:FDCalibrationOriginal} by the following gradient descent with diminishing stepsize algorithm
			\begin{itemize}
				\item Basic step: initialize $\rho^c_0, C_0$. A reasonable choice is $\rho^c_0 = 20$, which corresponds to the vehicles influencing one each other when the average distance among them is less than $50$ meters, and $C_0 = v^{\mathrm{limit}}_e\rho^c_0$, where $v^{\mathrm{limit}}_e$ is the speed limit on cell $e$ normalized by the sampling time $T$;
				\item $n$-th step: let $(\rho^c_n, C_n)$ descend along the gradient of the cost, namely
					\begin{align*}
						\rho^c_{n+1} 					& = \rho^c_{n} - \frac{\delta}{n}\nabla _{\rho^c}V_{(\rho^c, C)}\\
						C_{n+1} 						& = C_{n} - \frac{\delta}{n}\nabla_{C}V_{(\rho^c, C)}	
					\end{align*}
				with
					\begin{align*}
						\nabla_{\rho^c}V_{(\rho^c, C)} &=  \sum_{k\in\I^{\mathrm{FF}}(\rho^c_n)}\left(\varphi_k-\varphi_{(\rho^c_{n-1}, C_{n-1})}(\rho_k)\right)\frac{C_{n-1}}{(\rho^c_{n-1})^2}\rho_k\\ 
										&\quad				- \sum_{k\in\I^{\mathrm{C}}(\rho^c)}\left(\varphi_k-\varphi_{(\rho^c_{n-1}, C_{n-1})}(\rho_k)\right)C_{n-1}\frac{\rhojam-\rho_k}{(\rhojam-\rho^c_{n-1})^2}	\\
						\nabla_{C}V_{(\rho^c, C)}	  &= -\sum_{k\in\I^{\mathrm{FF}}(\rho^c_n)}\left(\varphi_k-\varphi_{(\rho^c_{n-1}, C_{n-1})}(\rho_k)\right)\frac{\rho_k}{\rho^c_{n-1}} \\ 
										&\quad	
						- \sum_{k\in\I^{\mathrm{C}}(\rho^c)}\left(\varphi_k-\varphi_{(\rho^c_{n-1}, C_{n-1})}(\rho_k)\right)\frac{\rhojam-\rho_k}{\rhojam - \rho^c_{n-1}}	\\	
							\I^{\mathrm{FF}}(\rho^c) &= \{k\in\K:\,0<\rho_k \leq \rho^c\}	\\
						\I^{\mathrm{C}}(\rho^c) &= \{k\in\K:\,\rho_k > \rho^c\}						
					\end{align*}
	where the gradients $\nabla_{\rho^c}V_{(\rho^c, C)}$ and $\nabla_{C}V_{(\rho^c, C)}$ are computed at $(\rho^c, C) = (\rho^c_{n-1}, C_{n-1})$, and $\delta >0$ is a fixed initial step size. Notice that if $\rho^c_n = 0$ then $\I^{\mathrm{FF}}(\rho^c) = \emptyset$, and conversely when $\rho^c_n = \rhojam$ then $\I^{\mathrm{C}}(\rho^c) = \emptyset$, and thus the previous summations are always well defined. Nonetheless, for numerical reasons, additional care should be taken in order to avoid $\rho^c< 0$ or $\rho^c>\rhojam$, for example projecting at each step $\rho^c_{n+1}$ into $[0, \rhojam]$ after the gradient update.
			\item Stopping criterion: stop if $\left\|\begin{bmatrix}\rho^c_n\\C_n\end{bmatrix}-\begin{bmatrix}\rho^c_{n-1}\\C_{n-1}\end{bmatrix}\right\| < \varepsilon$ for some small threshold $\varepsilon > 0$.
 			\end{itemize}	
			\item Calibration of the quadratic function in the congested region: the problem of calibrating the quadratic function for the congested region is cast into the quadratic problem: given the set of measures $\{(\rho_k^m, \varphi_k^m)\}_{k\in\K}$, $\K = \{1, \dots, K\}$, solve
			\begin{equation}
				\label{eq:FDCalibrationCongested}
				\begin{array}{ll}
				\min_{(a, b, c)}  	& 	\sum_{k\in\I^{\mathrm{C}}(\rho^c)}(\varphi_k - (a\rho_k^2+b\rho_k + c))^2\\
				\mathrm{s.t.}		&	a\left(\rho^c\right)^2 + b\rho^c + c_e = C	 \\
									& 	a\left(\rhojam\right)^2 + b\rhojam + c = 0 \\		
									&	a \geq 0
				\end{array}
			\end{equation}			
The problem \eqref{eq:FDCalibrationCongested} is computationally very simple and can be solved using off-the-shelf tools.
	\end{itemize}

Notice that as side products of the previous procedure we can compute the freeflow speed as $v_{ff} = C/\rho^c$ and, in case a bilinear Fundamental Diagram is also needed, the wave speed as $\omega = -C/(\rhojam-\rho^c)$, for each cell.

\subsubsection{Online density reconstruction algorithm}

We assume from now on that Fundamental Diagrams have either been calibrated or extended on the whole network, and that the matrix of splitting ratios has been pre-specified or estimated on the basis of field surveys.

We propose the following online algorithm for Density Reconstruction
	\begin{itemize}
		\item at the beginning of the $t$-th time slot, a centralized computation unit
			\begin{itemize}
				\item receives measurements $\{\varphi^m_e(t)\}_{e\in\E^m}$;
				\item \emph{flow estimation}: estimates the vector of outflows $\hatfout(t)$ by solving the following minimization problem
					\begin{equation}
						\label{eq:FindOutflows}
						\begin{array}{ll}
							\min_{\hatfout}  	& 	||(I-R^T)\hatfout||^2	+ \gamma\sum_{e\in\E^m}(\hatfout_e-\varphi^m_e(t))^2\\
							\mathrm{s.t.}		&	\hatfout \geq0
						\end{array}
					\end{equation}	
	Problem \eqref{eq:FindOutflows} aims to a) match outflows and measurements where available, by penalizing the squared difference between the two, and b) to balance outflows according to the splitting ratios. The latter term provides the estimate of flows on cells in which no measurement is available, and is performed \virg{as if} the network were at steady state, which is a simplifying assumption due to absence of a dynamical model for flows. The tunable parameter $\gamma$ selects whether more weight is given to matching estimated outflows and measurements (high $\gamma$), or to estimate the flows as if the network were at steady state (low $\gamma$). Once Problem \eqref{eq:FindOutflows} is solved, the vector of estimate of the inflows is easily computed according to Eq.~\ref{eq:fin}, by setting $\hatfin(t) = R^T\hatfout(t)$.
				\item receives the measurements $\{v^m_e(t)\}_{e\in\E}$ when available, or holds the last speed measurements received;
				\item For each cell $e$, computes the two possible densities $\rho_e^1$ (freeflow) and $\rho_e^2$ (congested) corresponding to flow $\hatfout_e(t)$ assuming that the local flow $\varphi_e = \varphi_e(\rho_e)$ is exactly determined by the Fundamental Diagram;
				\item For each cell $e$, computes the two velocities $v_e(\rho_e^1) = \frac{\hatfout_e}{\rho_e^1}$ (freeflow) and $v_e(\rho_e^2) = \frac{\hatfout_e}{\rho_e^2}$ (congested);
				\item Selects
					$$
						\tilde{\rho}_e(t) = \mathrm{arg}\min_{i=1,2}\{|v_e(\rho_e^i) - v^m_e(t)|\}
					$$
					as a rough estimate of the density. This estimate is only based on the actual measurements of flow and speed, and is generally very noisy, especially when the cell is in congestion. Therefore, the algorithm does not directly uses it;
				\item \emph{density estimation}: for each cell $e\in\E$, lets the density estimate evolve according to the observer equation \eqref{eq:observer}.
			\end{itemize}	
	\end{itemize}

Notice that in the proposed solution we avoid using Kalman filter based strategies as, due to the high degree of nonlinearity and uncertainty of the system. Analysis and minimization of the error variance of the proposed solution is left for future research.

\section{Optimal Sensor Placement}
\label{sec:optimalProblemFormulation}

The second problem we tackle in this paper is Optimal Sensor Placement, namely, the problem of deciding the position of sensors yielding a good trade off between performance and cost. Since assessing in a theoretical way the performance of algorithms for density reconstruction is difficult due to the nonlinear nature of the system, we simplify the setting and limit our attention to estimation of cumulative flows, namely, of the total outflows from the cells. The resulting, static, problem is then considered as a proxy for the more complicated problem built on the \emph{dynamic} density model.

We start by deriving some properties of cumulative outflows, and we proceed describing a simple linear model for cumulative flows estimation. This will help us formalizing the Optimal Placement problem.

\subsection{Flow linear constraints}
\label{subsec:flow_linear_constraints}

Let $f_e := \sum_{k=t_0}^{t_1-1}\fout_e(k)$ be the \emph{cumulative outflow} from cell $e$, namely, the total flow through the cell over the period of time $[t_0, t_1]$, and let $f = \begin{bmatrix}f_1&\dots&f_{|\E|}\end{bmatrix}^T\in\mathbb{R}_+^{\E}$ be the vector of cumulative outflows. By integrating the system's dynamics we have
	$$
		\ell_e(\rho_e(t_1) - \rho_e(t_0)) = \sum_{j\in\E}R_{je}f_j - f_e, \qquad e\in\E\setminus\R^i
	$$
Assume now that $[t_0, t_1]$ is a period of time whose duration is high enough, and that at both times $t_0$ and $t_1$ the number of vehicles in the network is low, for example, assume that $t_0$ and $t_1$ correspond to two consecutive midnights. Then the magnitude of the vector of differences of vehicles $\{\ell_e(\rho_e(t_1) - \rho_e(t_0)\}_{e\in\E\setminus\R^i}$ is small if compared
to the cumulative flows in the network, and the following relation holds approximately 
	\begin{equation}
		\label{eq:flowInKernel}
		\bar{L}f \approx 0
		\,,
	\end{equation}
where $\bar{L} \in \mathbb{R}^{\E\setminus\R^i\times \E}$, is the matrix obtained by removing from $L = R^T-I$ the rows corresponding to onramps. This imposes a linear constraint on the cumulative flows that we shall exploit in the next subsection.

\subsection{Linear measurement model and the Optimal Sensor Placement problem}
\label{subsec:measurementmodel}

In this subsection we study the performance of a linear estimator of the cumulative outflows and we show how to formalize the problem of Optimal Sensor Placement. 

Let $\E_m\subseteq\E$ be a set of cells in which sensors are placed. We assume the following simple linear measurement model
	\begin{equation}
		\label{eq:measurementOrig}
		y = H_{\E^m}f + \omega^f
	\end{equation}
where 
	\begin{itemize}
		\item $y_s$ is the measurement of the $s$-th sensor, namely $f_e + \eta_s$ if the $s$-th sensor is located on cell $e$;
		\item $H_{\E^m} \in \{0,+1\}^{p\times n}$, $p =  |E^m|$, $[H_{\E^m}]_{se} = 1$ if the $s$-th sensor is located on cell $e$, and $[H_{\E^m}]_{se} = 0$ otherwise, so that $H_{\E^m}\ones = 
		\ones$ and $\ones^TH_{\E_m}\ones = p$;
		\item $\omega^f$ is a random noise vector with zero mean and covariance matrix $\Sigma_{\mathrm{nom}}$, related to the measurement noise $\omega^\varphi$ described in the previous sections. For sake of simplicity, we shall often assume that the components of the noise, one for each sensor, are independent with same variance $\sigma_{\mathrm{nom}}^2$, so that $\Sigma_{\mathrm{nom}} = \sigma_{\mathrm{nom}}^2I$.
	\end{itemize}

Let now $n = |\E|$ and $r =\mathrm{rank}\{\bar{L}\}$, and consider a matrix $V \in \mathbb{R}^{n\times r}$ whose columns are an orthonormal basis of the right kernel of $\bar{L}^T$, i.e., $\bar{L}^TV = 0$ and $V^TV = I$. From \eqref{eq:flowInKernel} we get (approximatively) $f = Vz$ for some $z\in\mathbb{R}^r$, so that the measurement model can be rewritten as 
	\begin{equation}
		\label{eq:measurement}
		y = H_{\E^m} Vz + \omega^f
		\,.
	\end{equation}
	
Given $y$, consider a linear estimator of $z$, $\hat{z} = K_zy + q_z$, where $K_z\in\mathbb{R}^{r\times p}$ and $q_z\in\mathbb{R}^{r}$. The Best (minimum variance) Linear Unbiased Estimator of $z$ corresponds to the solution to
	\begin{equation}
	\label{eq:minimumVarianceProblem}
		\begin{array}{ll}	
			\min_{K_z,q_z}	 	& \mathbb{E}[(z-\hat{z})(z-\hat{z})^T]	\\
			\mathrm{s.t.}		& \mathbb{E}[z-\hat{z}] = 	0	\\
								& \hat{z} = K_zy + q_z
		\end{array}
	\end{equation}

The following Lemma formulate an equivalent problem. The proof is straightforward.

\begin{lemma}
\label{lemma:minimumVarianceProblem}
The solution $(K_z, q_z)$ to \eqref{eq:minimumVarianceProblem} is given by $q_z = 0$ and $K_z$ the solution of 
	\begin{equation}
	\label{eq:originalMinimumVarianceProblem}
		\begin{array}{ll}
			\min_{K_z}	 	& K_z\Sigma_{\mathrm{nom}}K_z'	\\
			\mathrm{s.t.}	& K_zH_{\E_m}V = I
		\end{array}
	\end{equation}
which is
	$$
		K_z = (V^TH_{\E_m}^T\Sigma_{\mathrm{nom}}^{-1}H_{\E_m}V)^{-1}V^TH_{\E_m}^T\Sigma_{\mathrm{nom}}^{-1}
		\,,
	$$
with error covariance
	$$
	\mathbb{E}[(z-\hat{z})(z-\hat{z})^T] = (V^TH_{\E_m}^T\Sigma_{\mathrm{nom}}^{-1}H_{\E_m}V)^{-1}	
	$$
\end{lemma}	

An immediate consequence of Lemma~\ref{lemma:minimumVarianceProblem} is that the Best Linear Unbiased Estimator (BLUE) of $f$ is
	$$
		\hat{f} = K_fy = V(V^TH_{\E_m}^T\Sigma_{\mathrm{nom}}^{-1}H_{\E_m}V)^{-1}V^TH_{\E_m}^T\Sigma_{\mathrm{nom}}^{-1}y
	$$
and its error covariance is 
	\begin{align*}
		V_p(\E^m) 
			&	=	\mathbb{E}[(f-\hat{f})(f-\hat{f})^T]	\\
			& 	= V(V^TH_{\E^m}^T\Sigma_{\mathrm{nom}}^{-1}H_{\E^m}V)^{-1}V^T
		\,.
	\end{align*}

The quantity $V_p(\E^m)$ depends on a) the (right kernel of the) matrix of splitting ratios via the matrix $V$ and the nominal variance of the noise $\omega^f$, two parameters that are assigned, and b) on the locations of the sensors, which is the set $\E^m$, via the matrix $H_{\E^m}$. For this reason, we will take the magnitude of $V_p(\E^m)$, measured via its trace, as our metric to measure the performance of a sensor network placed on the cells $\E^m$. 

Clearly, with no additional constraints the optimal placement is to equip every cell with sensors. This is straightforward as equipping all cells means setting $H_\E = I$, and from $H_\E^TH_\E = I \geq H_{\E^m}^TH_{\E^m}$ immediately descends $V_p(\E) \leq V_p(\E^m)$, for any $\E^m\subseteq\E$.

Each device has however a non-negligible purchase and maintenance cost, which has to be considered when designing the sensor network. For sake of simplicity, in this paper we make the simplifying assumption that the cost of a network over its lifetime is proportional to its number of sensors via a coefficient $c > 0$, so that the cost of deploying sensors on $\E^m$ is $c|\E^m|$.

We thus consider the following

\begin{problem}[Optimal Sensor Placement]
Let $\G = (\V, \E)$ be a traffic network with splitting ratios $R$ and cumulative flows noise variance  $\sigma_{\mathrm{nom}}^2$. Find $\E^m$ which solves
	\begin{equation}
	\label{eq:originalCombinatorialProblem}
		\begin{array}{ll}	
			\min_{\E^m}	 	& \trace{V_p(\E^m)} + c|\E^m|
		\end{array}		
	\end{equation}
\end{problem}

The optimal $\hat{\E}^m$, solution to \eqref{eq:originalCombinatorialProblem}, trades off between the network performance, which is measured by the trace of the estimator error covariance, and the total cost of the network. Clearly, the two have a contrasting effect on the number of deployed sensors. It is however inherently combinatorial, the optimal position of the sensor being in general hard to find and requiring an exhaustive search among all the possibilities, which is intractable even for relatively small network dimensions. 

We approach the problem by proposing an heuristic that relaxes it into a convex problem. Such a strategy is described after the following brief discussion on the minimum required number of sensors.

\subsection{Minimum number of sensors}
\label{subsec:minimumNumber}

Before proposing our method for solving \eqref{eq:originalCombinatorialProblem}, we prove that there exists a lower bound on the number of sensors $|E^m|$ in order the $\trace{V_p(\E^m)}$ to be finite. As it will be proven, below such number, which corresponds to the number of onramps of the system, the problem of reconstruction of flows admits infinite solutions. 

To this aim, relabel the cells in such a way that onramps are the first $1, \dots, |\R^i|$ cells so that we can partition $\bar{L}$ as
	$$
		\bar{L} = \begin{bmatrix}L_{on}&L_{nn}\end{bmatrix}
	$$
where $L_{nn}$ model the mutual influence of flows on non-onramp cells, and $L_{on}$ models the influence of onramps on non-onramp cells.

Define now a dual graph $\mathcal{G}^d = (\V^d, \E^d)$ in which the roles of cells and junctions are in a way reversed, and in particular in which $\V^d = \E\setminus \R^i$ and $(e, j) \in \E^d$ if $e\neq j$ and $[L_{nn}]_{ej}\neq 0$. Then it is easy to see that $L_{nn}^T$ is a sublaplacian of $\G^d$. The following result is adapted from \cite{LovisariCDC:14}.

\begin{lemma}
\label{lemma:technical}
Let $\mathcal{G} = (\V, \E)$ be a graph and $J\in\mathbb{R}^{|\V|\times|\V|}$ be a sublaplacian of $\G$. Then all the eigenvalues of $J$ have negative real part except possibly eigenvalues in $0$. Moreover, if $\S$ is the set of nodes $v$ for which $\sum_{u}J_{vu} < 0$, then $J$ is invertible if for every $u$ there exists a directed path in $\G$ from $u$ to a node $v\in\S$.
\end{lemma}

In the case under analysis, we take $J = L_{nn}^T$ and $\S$ to be the set of cells directly following an onramp, namely $\S = \{e\in\E:\,\exists j\in\R^i:\,R_{je} > 0\}$. Recall that by assumption for every cell $e\not\in\R^i$ there exists an origin $j\in\R^i$ and a path from $j$ to $e$, so there must also be a $k\in\S$ and a path from $k$ to $e$ (at most being $k=e$), so the assumptions of Lemma~\ref{lemma:technical} are satisfied. This establishes that $L_{nn}$ is invertible, and therefore that $\bar{L}$ is a full row-rank matrix with rank $|\E\setminus\R^i|$. Since the number of columns of $\bar{L}$ is $|\E|$, it follows that its kernel has dimension $r = |\R^i|$, and thus $\rank{V} = |\R^i|$. 

We offer two interpretations on this fact:
	\begin{itemize}
		\item Ideal measurement scenario: assume $p = |\E^m|$ ideal measurements at sensors $s(1), \dots, s(p)$ are available, $y_{s(i)} = f_{s(i)}$, $i=1, \dots, p$. Then a solution to the system of equations
			$$
				\begin{cases}
					\bar{L}\hat{f} = 0	\\
					\hat{f}_{s(i)} = f_{s(i)}, i = 1, \dots, p
				\end{cases}
			$$
is a candidate vector of cumulative flows. Then, if $p\geq |\R^i|$ the system has a unique solution which is the true vector of flows. Conversely, if $p < |\R^i|$, the system is undetermined;
		\item Noisy measurement scenario: assume $p = |\E^m|$ noisy measurements at sensors $s(1), \dots, s(p)$ are available and assume to adopt the Best Linear Unbiased Estimator presented above to estimate the flows. If $p = |\E^m| < |\R^i|$, then $\rank{H_{\E^m}^TH_{\E^m}} < |\R^i| = \rank{V}$, which implies that $\rank{V^TH_{\E^m}^TH_{\E^m}V} < |\R^i|$. However, $V^TH_{\E^m}^TH_{\E^m}V \in \mathbb{R}^{|\R^i|\times |\R^i|}$, so the matrix is singular, and therefore the trace of the error covariance is unbounded. Conversely, if $p \geq |\R^i|$  then $\rank{V^TH_{\E^m}^TH_{\E^m}V} = |\R^i|$ and the trace of the error covariance is bounded.
	\end{itemize}

\subsection{Relaxation via Virtual Variances}
\label{subsec:solution}

The solution that we propose is based on the observation that cells that are not endowed with sensors can be interpreted as cells in which sensors have infinite noise variance. It turns out that an equivalent formulation of \eqref{eq:originalCombinatorialProblem} is
	\begin{equation}
	\label{eq:slightlyChangedCombinatorialProblem}
		\begin{array}{ll}	
			\min_{\E^m}	 	& \trace{V(V^T\Sigma V)^{-1}V^T} + c\sigma_{\mathrm{nom}}^2\ones^T\Sigma^{-1}\ones\\
			\mathrm{s.t.}	& \Sigma_{ee} = 
									\begin{cases}
										+\infty, &e\not\in\E^m	\\
										\sigma_{\mathrm{nom}}^2, &e\in\E^m
									\end{cases}
		\end{array}		
	\end{equation}
In fact
	\begin{itemize}
		\item the term $H_{\E^m}$, which represents which cells are endowed with a sensor, is the identity matrix, namely, all cells are endowed with a sensor -- except, some of them have infinite noise variance and thus provide no information;
		\item the second term in the cost corresponds to $c|\E^m|$ as
			$$
				c\sigma_{\mathrm{nom}}^2\ones^T\Sigma^{-1}\ones
				=
				c\sigma_{\mathrm{nom}}^2\left(\sum_{e\in\E^m}\frac{1}{\sigma_{\mathrm{nom}}^2}\right)
				=
				c\sum_{e\in\E^m}1 = c|\E^m|
			$$
	\end{itemize}

In other terms, \eqref{eq:slightlyChangedCombinatorialProblem} corresponds to assigning a \emph{virtual variance} $\Sigma_{ee} = \sigma_e^2$ to each sensor, and decide for which it should be $\sigma_e^2 = \sigma_{\mathrm{nom}}^2$, the cells in $\E^m$, and for which it should be $\sigma_e^2 = + \infty$. 

We call $\Sigma$ the (diagonal) matrix of virtual variances, and we let the corresponding trace of error covariance be denoted, with an abuse of notation, $V_p(\Sigma) = \trace{V(V^T\Sigma V)^{-1}V^T}$. 
\medskip

Our approach is then based on the intuitive idea that increasing the variance on the sensors that are not crucial for the solution of \eqref{eq:originalCombinatorialProblem} should not have a strong effect on the performance term $V_p(\Sigma)$. More formally, we consider the following relaxed version of the previous problem
	\begin{equation}
	\label{eq:firstRelaxedProblem}
		\begin{array}{ll}	
			\min_{\Sigma\in\mathbb{D}_n}	 	& \trace{V(V^T\Sigma^{-1}V)^{-1}V^T} + c\sigma_{\mathrm{nom}}^2\ones^T\Sigma^{-1}\ones\\
			\mathrm{s.t.}						& \Sigma_{ee} \geq \sigma_{\mathrm{nom}}^2, \forall e\in\E
		\end{array}		
	\end{equation}
where $\mathbb{D}_n$ is the set of diagonal matrices of dimension $n$.

We now slightly rewrite the cost. First, by the well known property of trace $\trace{AB} = \trace{BA}$ and $V^TV = I$, it follows $\trace{V(V^T\Sigma^{-1}V)^{-1}V^T} =  \trace{(V^T\Sigma^{-1}V)^{-1}}$. Second, we make the change of variables $\Omega = \Sigma^{-1}$. In this way, we obtain the following problem
	\begin{equation}
	\label{eq:secondRelaxedProblem}
		\begin{array}{ll}	
			\min_{\Omega\in\mathbb{D}_n}	 	& \trace{(V^T\Omega V)^{-1}} + \gamma\ones^T\Omega\ones\\
			\mathrm{s.t.}		& 0\leq \Omega \leq \Sigma_{\mathrm{nom}}^{-1}
		\end{array}		
	\end{equation}
where $\gamma$ is a tunable parameter. The choice $\gamma = c\sigma_{\mathrm{nom}}^2$ is a natural one due to the previous discussion, but since $\gamma$ influences the relative weight of performance (penalized for low $\gamma$) and cost (penalized for high $\gamma$) we leave it as an additional degree of freedom. Notice, finally, that the term $\ones^T\Omega\ones$ corresponds to the $\ell_1$ norm of the inverse of the variances, a term which is commonly used term to sparsify solutions of optimization problems.

The \emph{Virtual Variance algorithm} proceeds as follows:
	\begin{enumerate}
		\item solve \eqref{eq:secondRelaxedProblem} and denote by $\Omega$ its solution;
		\item compute $\Sigma^{-1} = \Omega$
		\item discard all cells whose virtual variance is above a fixed discard threshold $\mathcal{T}_d$
	\end{enumerate}
 	
As explained above, if the found solution provides high virtual variances at locations where sensors are redundant then this is effectively a way to select the most important cells where to place sensors.

However, this strategy does not, in general, provide good solutions to the problem. Indeed, numerical simulations have shown that in the considered scenario the solution of \eqref{eq:secondRelaxedProblem} can be often interpreted as endowing all cells of the network with sensors with average virtual variance, rather than keeping it low in some of them and high in others.

In order to enhance \emph{diversity} between sensors, we enrich the cost of \eqref{eq:secondRelaxedProblem} with a term that aims to penalize homogeneity. This is reminiscent of \emph{dissensus} (as opposed to consensus) strategies in multi agent networks, in which each agent possesses a value and the goal of the network is to differentiate as much as possible such values.

In this paper, we make the following simple choice. Let $W \in \mathbb{R}^{n\times n-1}$ be an orthonormal base of the subspace orthogonal to $\ones$, that is, $W^*\ones = 0$ and $W^*W = I$. We add to the cost in \eqref{eq:secondRelaxedProblem} a term that is proportional to $e^{-\ones^TW^*\Omega\ones}$. Since the columns of $W$ span the orthogonal to $\ones$, $W^*\Omega\ones$ is high when the element on the diagonal of $\Omega$, which are gathered in the vector $\Omega\ones$, are different one with respect to the other, and is low otherwise.

We propose the following optimization problem
	\begin{equation}
	\label{eq:relaxedProblem}
		\begin{array}{ll}	
			\min_{\Omega\in\mathbb{D}_n}	& \trace{(V^T\Omega V)^{-1}} + \gamma\ones^T\Omega\ones + \kappa e^{-\ones^TW^*\Omega\ones}\\
			\mathrm{s.t.}					& 0\leq \Omega \leq \sigma_{\mathrm{nom}}^{-2}I
		\end{array}		
	\end{equation}
which, notice, is \emph{convex} in the diagonal entries of $\Omega$. Here $\gamma$, the \emph{total variance weight}, and $\kappa$, the \emph{discrepancy weight}, are tunable parameters. Notice that high $\gamma$ penalizes the number of sensors, thus yielding to solutions with higher virtual variances at the expense of poor performance. 

We shall provide in the following paragraphs examples of application of the previous optimization problem for different sensor placement scenarios. We anticipate here that one can observe, as required, that the virtual variances solution of \eqref{eq:relaxedProblem} are distributed in a strongly bimodal way, with low and high values being different by several orders of magnitude. As a consequence, it is usually easy to distinguish among the two groups and discard cells whose contribution to the performance metric would be negligible.

\subsubsection{Optimal Sensor Placement with geographical or budget constraints}

In this paragraph we discuss two variations of the previous procedure, which address the additional problems of geographical constraints and of strict budget limitation.
\medskip 

\emph{Optimal Sensor Placement with geographical constraints}
\medskip 

The first scenario we address is concerned with the scenario in which in which a) some cells cannot be endowed with sensors, for example for physical reasons, and/or b) subsets of cells for which either all cells are endowed with sensors, or none are. An example of the latter constraints is a multi-lane road which is modelled using parallel cells (on different lanes) and on which the traffic manager can deploy induction loop. The former are buried underground are are usually required to cover the whole carriageway. As such, if $i$ and $j$ are, for example, two parallel cells on the two lanes on a certain section or road, then either $i$ and $j$ are both equipped with a sensor (i.e., the induction loop), or not.

To encompass this type of constraints in the proposed procedure, let $\E_{am} \subseteq \E$, $|\E_{am}| = \kappa$, be the subset of available cells, and let $H_{\E_{am}} \in\{0,+1\}^{\kappa\times n}$ be built as in Subsection~\ref{subsec:measurementmodel}. Further, let $\S \subseteq \mathcal{P}^{\E}$ be the set of subsets of $\E$ which must be simultaneously equipped, or not, with sensors, where $\mathcal{P}^{\E}$ is the powerset, or set of subsets, of $\E$.

Then the following problem
	\begin{equation}
	\label{eq:relaxedProblemConstrained}
		\begin{array}{ll}	
			\min_{\Omega\in\mathbb{D}_k}	& \trace{(V^TH_{\E_{am}}^T\Omega H_{\E_{am}}V)^{-1}} + \gamma\ones^T\Omega\ones + \kappa e^{-\ones^TW^*\Omega\ones}\\
			\mathrm{s.t.}					& 0\leq \Omega \leq \sigma_{\mathrm{nom}}^{-2}I			\\
											& \Omega_{ii} = \Omega_{jj}, \forall i,j\in \sigma, \forall \sigma\in\S
		\end{array}		
	\end{equation}
is \eqref{eq:relaxedProblem} once we constrain sensors to be place on cells in $\E_{am}$ only and simultaneous sensor placement to happen according to the constraints specified by the set $\S$. Notice that the diagonal entries of the solution $\Omega$ are the inverse of the virtual variances on the cells in $\E_{am}$ only. The matrix $W$ is defined as previously, but with suitable dimension ($\kappa\times \kappa-1$). As in the general problem, cells are chosen only if the corresponding virtual variance is below a certain threshold, and clearly cells that are not in $\E_{am}$ cannot be chosen. Notice that the constraint $ \Omega_{ii} = \Omega_{jj}$ in \eqref{eq:relaxedProblemConstrained}, or the less requiring $|\Omega_{ii} - \Omega_{jj}| \leq \varepsilon$, for some tunable parameter $\varepsilon$, are convex constraints, so \eqref{eq:relaxedProblemConstrained} remains convex.

\begin{remark}
\label{rem:MinEam}
By the discussion in Subsection~\ref{subsec:minimumNumber}, the minimum number of sensors is $r = |\R^i|$. As such, if $|\E_{em}|< |\R^i|$ the problem \eqref{eq:relaxedProblemConstrained} is not well posed and the solution will only have very high virtual variances. Clearly, such a solution is not acceptable and should be discarded.
\end{remark}

\subsubsection{Optimal Sensor Placement with budget constraints}

In this second scenario we discuss budget constraints in the form of constraints on the maximum number of chosen sensors, a very common requirement in real-case applications.

We propose a solution based on the following iterative approach: 
	\begin{itemize}
		\item Initialization: set $\gamma(0)$ and $\kappa$ to some prespecified nonnegative values; $t_{\mathrm{max}}$ to the maximum number of iterations; $n_{\mathrm{max}}$ to the maximum number of sensors;
		\item $t$-th step
			\begin{itemize}		
				\item The problem \eqref{eq:relaxedProblem} is solved with $\gamma = \gamma(t)$;
				\item Let $n(t)$ be the number of sensors in the solution provided by the Virtual Variance algorithm. Then
					\begin{itemize}
						\item If $n(t)\leq n_{\mathrm{max}}$, or if $t \geq t_{\mathrm{max}}$, the procedure stops;
						\item Otherwise, set $\gamma(t+1) = g(\gamma(t))$, where $g$ is a monotonically increasing increasing function, and the procedure iterates.
					\end{itemize}
		\end{itemize}
	\end{itemize}

The rationale behind this procedure is that, as previously discussed, $\gamma$ penalizes a low total sum of the virtual variances. Therefore, by iteratively increasing $\gamma$ the solution to \eqref{eq:relaxedProblem} will tend to exhibit more and more high virtual variances, thus reducing the number of sensors.

\begin{remark}
\label{rem:MinNumSen}
Once again, and related to the discussion in Remark~\ref{rem:MinEam}, the specified maximum number of sensors cannot be be less than $r = |\R^i|$ by the results presented in Section~\ref{subsec:minimumNumber}. If this is not the case, numerical experiments show that the algorithm either simply iterates until the number of iterations reaches $t_{\mathrm{max}}$, or the found solution exhibits an extremely high $\trace{(V^TH_{\E_{am}}^T\Omega H_{\E_{am}}V)^{-1}}$ - due to the fact that the internal matrix is (numerically) almost singular. As in Remark~\ref{rem:MinEam}, such a solution is not acceptable and should be discarded.
\end{remark}

\section{Numerical experiments}
\label{sec:numericalExperiments}

\subsection{Numerical Experiments for the Optimal Sensor Placement}
\label{subsec:numerical_osp}

In this subsection we present the results of two numerical experiments. In the first, we solve the problem of Optimal Sensor Placement in a small (but not trivial) regular grid, in the second, we apply the procedure to the real-world case of the Rocade Sud. 

\subsubsection{Regular grid}

\begin{figure}
\centering 
\includegraphics[width=0.4\textwidth]{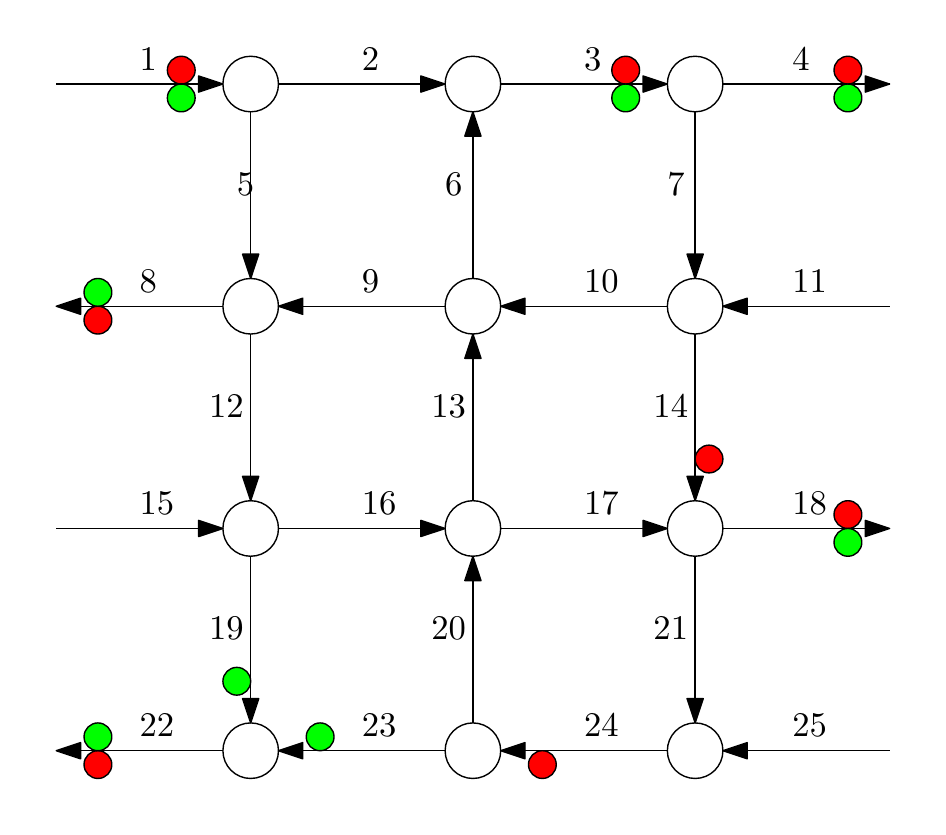}
\caption{The regular grid network used in the numerical experiment. The green dots correspond to the $8$ cells selected by the Virtual Variance algorithm with $\gamma = 2$ and $\kappa = 20$. The red dots correspond to cells selected via exhaustive search when the number of possible sensors is $8$.}
\label{figure:ExampleGrid}
\end{figure}

Consider the regular grid composed of $25$ cells shown in Figure~\ref{figure:ExampleGrid}. We want to solve the problem \eqref{eq:originalCombinatorialProblem} by \eqref{eq:secondRelaxedProblem} with parameters $\sigma_{\mathrm{nom}}^2 = 1$, nominal sensor variance, and $c=1$, cost of a single sensor.

The network is small enough to run an exhaustive search to solve the problem \eqref{eq:originalCombinatorialProblem}. In particular, for each $h = 4, 5, 6, \dots, 21$, we try all possible combinations of $h = |\E^m|$ sensors, thus finding the one that minimizes $V_p(\E^m)$. 

We also run our Virtual Variance algorithm with total variance weight $\gamma = 2$ and discrepancy weight $\kappa = 20$. Further, we set the threshold on the virtual variance for discarding a sensor to $\mathcal{T}_d = 100$.

\begin{figure}
\centering 
\includegraphics[width=0.46\textwidth]{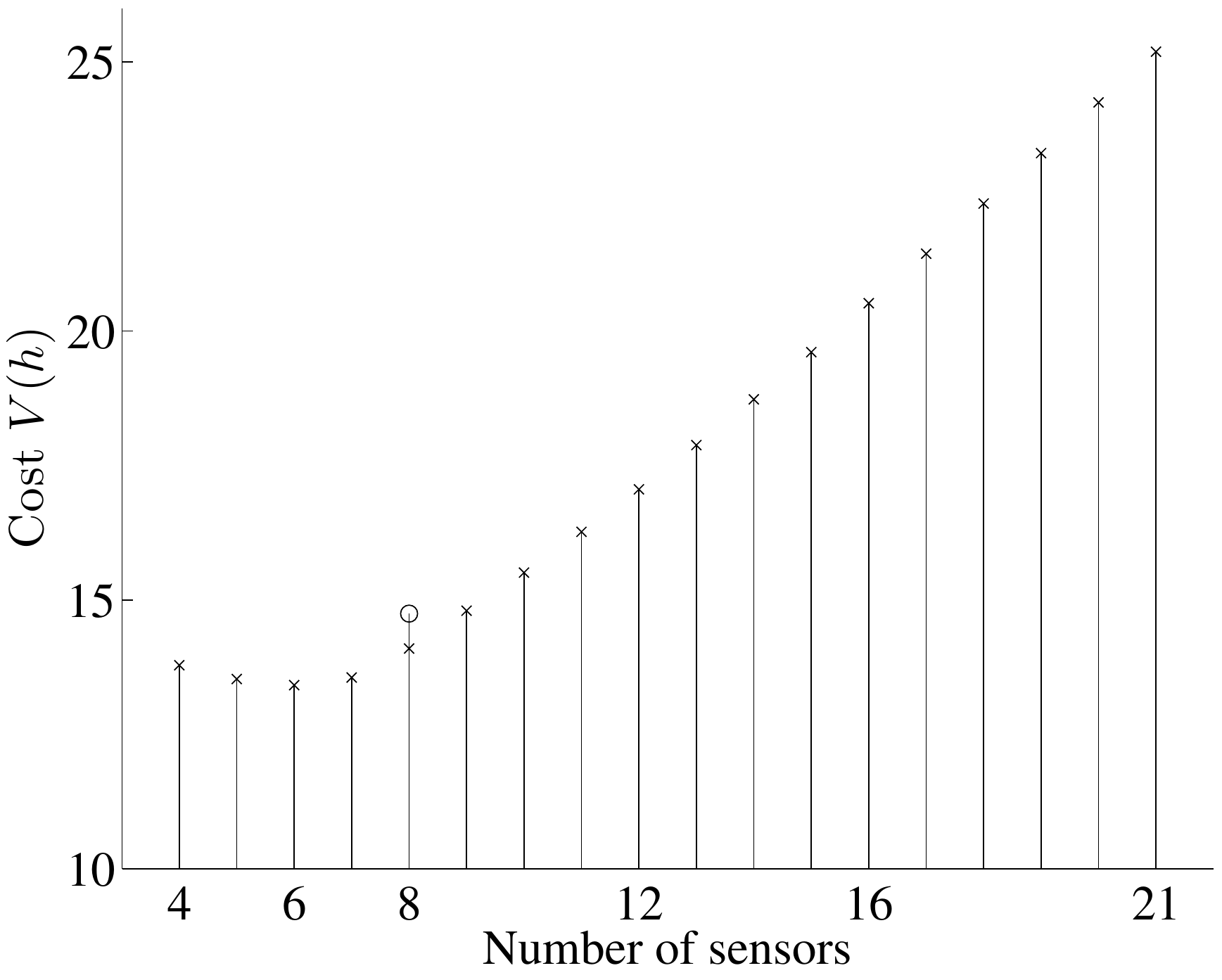}
\caption{Results of the exhaustive search and of the virtual variance algorithm.}
\label{figure:ResultsExtensiveVSVirtualVariance}
\end{figure}

We illustrate the results in Figures~\ref{figure:ExampleGrid} and \ref{figure:ResultsExtensiveVSVirtualVariance}. In Figures~\ref{figure:ExampleGrid} the $8$ cells chosen by the Virtual Variance algorithm are represented as a green dot marks, and the best possible placement with $8$ sensors (found by exhaustive search) as red dot marks. It can be noticed that both procedures place most of sensors at the boundaries of the network.

Figure~\ref{figure:ResultsExtensiveVSVirtualVariance} shows instead the total cost of the best placement obtained through exhaustive search for $h = 4, \dots, 21$, and the total cost found by the Virtual Variance algorithm with $8$ sensors. By total cost we mean the sum of estimation performance and network cost $V_p(\E^m) + c|\E^m|$, which, notice, is \emph{not} the metric that is used in the Virtual Variance algorithm. Nonetheless, it is appreciable that the Virtual Variance algorithm not only provides a solution whose number of sensors is close to the global optimum (which is with $6$ sensors), but also that, using $8$ sensors, the Virtual Variance algorithm places them almost in the optimal way.

\subsubsection{Rocade Sud}
\label{subsubsec:RocadeSud}

\begin{figure}
\centering 
\begin{tabular}{c}
\includegraphics[width=0.45\textwidth]{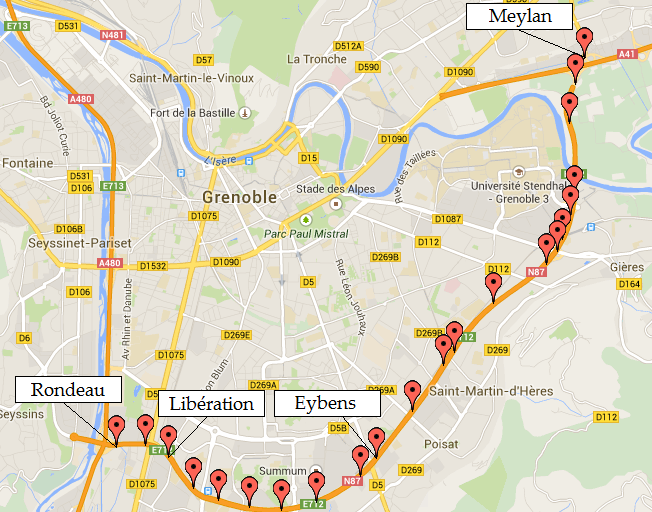}~
\includegraphics[width=0.5\textwidth]{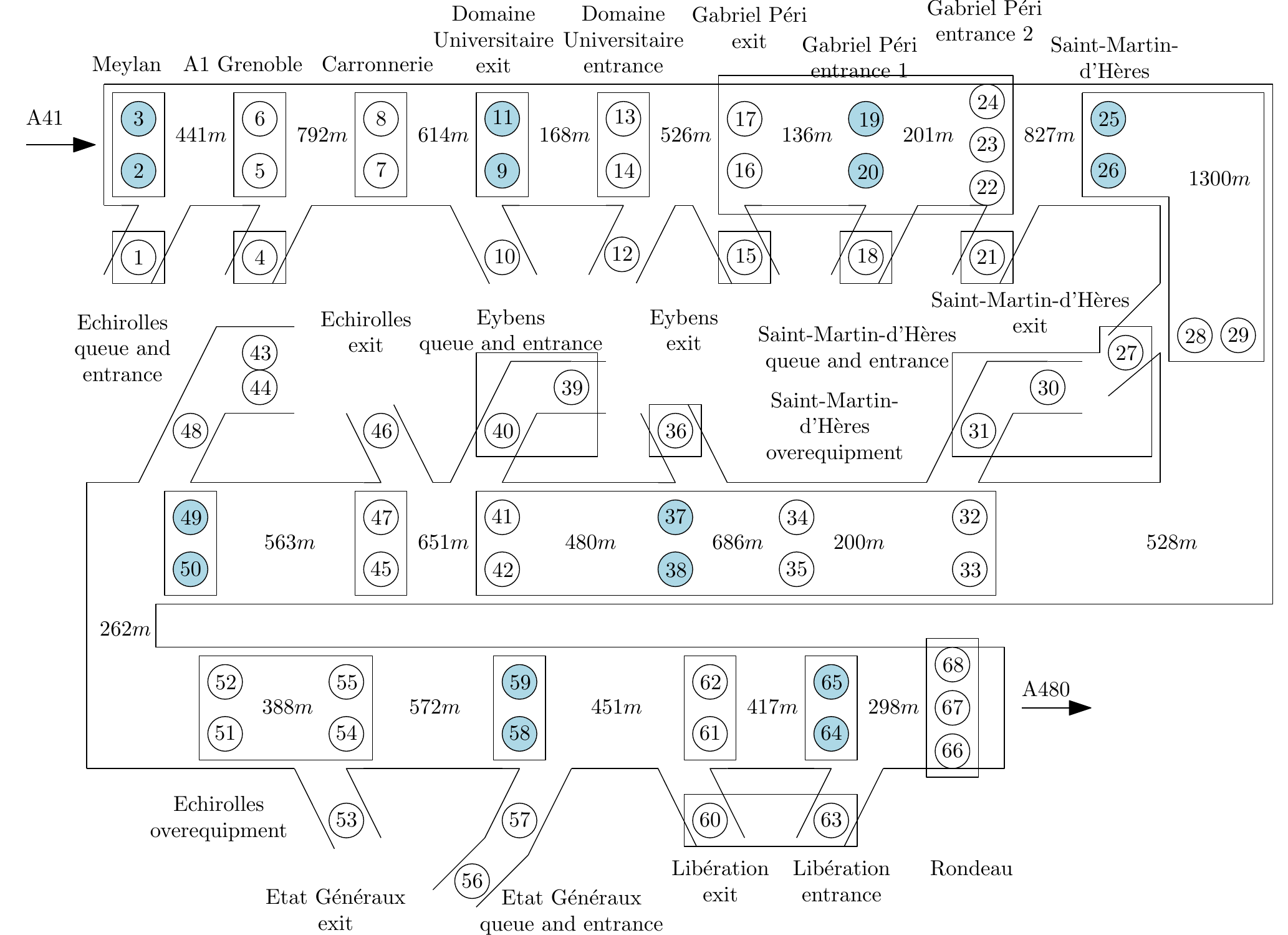}
\end{tabular}
\caption{The experimental setting: the town of Grenoble and the Rocade Sud and a stylized version of the freeway. The positions along the main line of the $22$ sections of the main line in which sensors have been placed is shown as red pin.
The positions of the $68$ fixed sensing locations are shown in the stylized map. Each sensing location also corresponds to a cell. Light blue circles denote fixed sensors that are selected by the Virtual Variance algorithm and are used in the implementation of the Density Reconstruction algorithm. Each rectangle represents one FCD segment, often providing average speed measurements over more than one cell.}
\label{figure:RocadeSudStylized}
\end{figure}

Our second experimental setting is the Grenoble Traffic Lab (GTL), a network of sensors deployed for monitoring and research purposes along the Rocade Sud, a peri-urban $10.5$ km long freeway connecting the two highways A41 (north-west) to A480 (south) in the town of Grenoble in the south of France, see Figure~\ref{figure:RocadeSudStylized}. The network is composed of $135$ magnetometers buried in the ground on both lanes along the main line every $250$ meters (on average), on each onramp and offramp, and on three connectors from urban roads to three onramps, for a total of $68$ \emph{sensing locations}. For our purposes, each sensing location will correspond to one pair of sensors. For a detailed report on the GTL, we refer to \cite{CanudasIEEECS:15}. Furthermore, and for sake of simplicity, while in the real network sensors are deployed in pairs, we shall assume from now on that each sensing location has only one sensor (as we shall always discard sensor speed measurements).

Figure~\ref{figure:RocadeSudStylized} shows the position of each of the $22$ sections of the main line in which there are sensing locations on both slow and fast lanes (and usually a ramp). In the same figure we also show a stylized representation of the freeway, including ramps and queues, the positions of the $68$ sensing locations, and the distance between consecutive measurement sections along the main line.

We partition the Rocade in cells in such a way each cell includes one sensor, so that in Figure~\ref{figure:RocadeSudStylized} each numbered circle also corresponds to one cell.

Here, we do not consider onramps and offramp, limiting our attention to the main line of the Rocade Sud. The reason for this choice is that the Rocade has 10 onramps along its main line, which, summed to the two cells in the very first section of road, imply a minimum number of sensor of $12$ by the discussion in Subsection~\ref{subsec:minimumNumber}. While this number is not high on its own, numerical experiments not reported in this paper have shown that good estimation performance require a number of sensors which is too high in most realistic (i.e., non academical) implementations. 

\begin{figure}
\centering %
\includegraphics[width=0.46\textwidth]{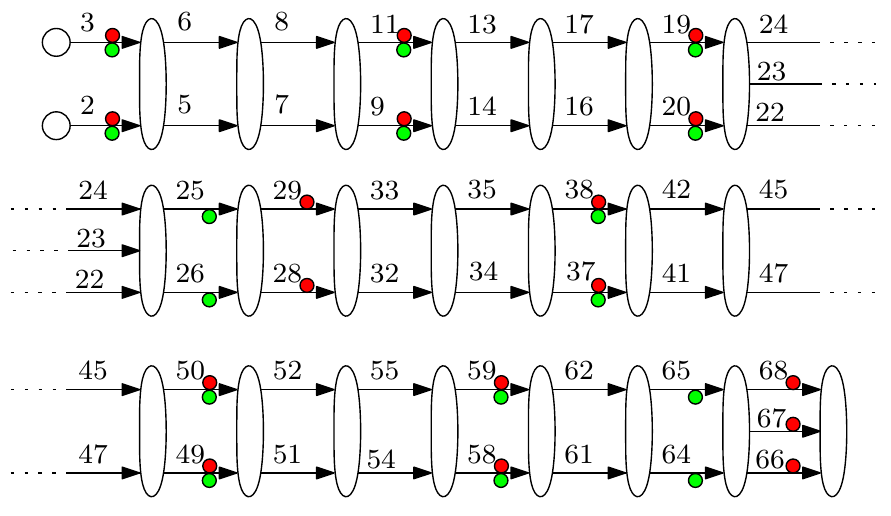}
\caption{Stylized representation of the main line part of the Rocade Sud. White ovals represents junctions of the graph. The selected positions of $17$ fixed sensors by CNIR are marked by red dots, those of the $16$ sensors chosen by the Virtual Variance algorithm by green dots.}
\label{figure:ExampleRocadeMainLinePositions}
\end{figure}

The corresponding reduced graph representation, essentially made of several groups of parallel edges, consists of $46$ cells. We provide a stylized version of it in Figure~\ref{figure:ExampleRocadeMainLinePositions}.

The matrix of splitting ratios could be estimated from the data, but we make here a simpler assumption and assume that vehicles split according to the following rule: vehicles on slow lane cells remain on slow lane or turn into fast lane with a $70$\%-$30$\% rule, and analogously for vehicles on the fast lane. If the next section has three parallel cells (cells 22-23-24 and cells 66-67-68), vehicles spit uniformly in such three cells.

We did not run an exhaustive search due to the relatively high dimension of the network and the consequent relevant computational load. Our comparison is instead with the locations of fixed loops installed by the Government Agency \emph{Centre national d'information routi\'{e}re (CNIR)} \cite{CNIR}, which correspond to cells marked with a red dot in Figure~\ref{figure:ExampleRocadeMainLinePositions}. We shall show that the latter are positioned in a way that is in good accordance with the results of the Virtual Variance algorithm, even though our procedure allows for a slightly better design.

We run the Virtual Variance algorithm in three scenarios: 1) unconstrained scenario with total variance weight $\gamma = 0.2$ and discrepancy weight $\kappa = 20$; 2) unconstrained scenario with $\gamma = 1$ and $\kappa = 20$; 3) constrained scenario with number of sensors at most $10$, initial $\gamma = 0.2$, and $\kappa = 20$. We assume that $\sigma_{\mathrm{nom}}^2 = 1$ and that the cost per sensor is $c = 1$. In all cases, sensors are constrained either to be present in both lanes on a same section, or to be absent.

The results are summarized in Figure~\ref{figure:ExampleRocadeMainLinePositions} and Table~\ref{table:resultsRocade}. The tables shows the number of sensors in the solution computed via the Virtual Variance algorithm, as well as the trace of the corresponding estimation error covariance $V_p(\E^m)$ and the total cost $V(\E^m) = V_p(\E^m) + c|\E^m|$. In Figure~\ref{figure:ExampleRocadeMainLinePositions}, cells found in the unconstrained scenario with $\gamma=0.2$ are marked with a green dot. As can be seen in Table~\ref{table:resultsRocade}, our algorithm requires one less sensors than the network deployed by CNIR and in addition the trace of the error covariance is smaller. 

In the constrained scenario and in the unconstrained scenario with high $\gamma$ (which, as explained above, indirectly penalizes the number of sensors), the trace of the error covariance $V_p(\E^m)$ increases, as expected. Furthermore, the chosen cells in the latter two cases are subsets of the cells chosen in the unconstrained case: in particular, in the constrained scenario all cells are kept except $13$, $14$, $49$ and $50$, and in the unconstrained scenario with high $\gamma$ the algorithm further discards cells $64$ and $65$ 

\begin{table}
\begin{center}
\begin{tabular}{l|cccc}
Scenario 			&	$\gamma$	& \begin{tabular}{c}optimal \\\# sensors\end{tabular}		& 	$V_p(\E^m)$ 		& 	$V(\E^m)$ 		\\
\hline
Fix							&				&	17			& 	$3.8161$	& 	$20.8161$	\\		
Unconstrained 				&	$0.2$		& 	14			&	$3.6867$	&	$17.6867$	\\ 
Unconstrained				&	$1$			&	10 			& 	$5.1732$	&	$15.1732$	\\
Constrained, \# $\leq 12$	&	$0.32$		&	12			& 	$4.4972$	&	$16.4972$	
\end{tabular}
\end{center}
\caption{Results of the four considered scenarios.}
\label{table:resultsRocade}
\end{table}

\subsection{Density reconstruction - experimental results}
\label{subsec:numericalExperimentsDensityReconstruction}

We provide here numerical results of the implementation of the data fusion algorithm for density reconstruction on data from the Rocade Sud.

On each sensing location and every $T = 15$ seconds, the system counts the number $\varphi^m_e$ of vehicles that crossed the location, their average speed $v^m_e$, and the average occupancy $o^m_e$ of the location. Since the latter is approximatively proportional to the density of vehicles, so we shall assume that sensors can directly measure densities. 
In addition to fixed sensors, we use Floating Car Data provided by INRIX. Following INRIX schema, the Rocade has been further partitioned into \emph{FCD segments}. One measurement of average speed is available on each FCD segment every $1$ minutes. FCD segments partition the whole main line of the Rocade and include most onramps and offramps, but single lanes are not distinguished along the main line. FCD segments are represented in Figure~\ref{figure:RocadeSudStylized} as rectangles encircling several sensing locations/cells.
For our experiments, we employ the sensor configuration obtained in the previous section via the Virtual Variance algorithm with $\gamma=0.2$. In particular, and in order to prove that the method shows good performance even with sparse equipment, we only use the sensors on the main line, which are, for reference, shown in light blue in Figure~\ref{figure:RocadeSudStylized}. Further, we don't use any information on flow or speed on the ramps.

Calibration of the Fundamental Diagram was performed via the algorithm described in Paragraph~\ref{subsubsec:offline_calibration} and using the data from the GTL sensor network from April 10th, 2014, a working day (a Thursday) exhibiting very standard traffic pattern: 
	\begin{itemize}
		\item very limited night time traffic;
		\item a peak of congestion in the morning (8:00 - 10:00), triggered by vehicles exiting towards the city from the Rocade at the offramp of Eybens (cells 37/38) and spilling back until Meylan (cells 2/3), and a second, smaller peak of congestion triggered by vehicles entering in A480 at Rondeau (cells 66/67/68) but blocked by the high traffic on A480, and spilling back until around Lib\'{e}ration (cells 61/62);
		\item a third, smaller, congestion triggered around Eybens around 14:00-15:00;
		\item in general, medium/heavy but fluid traffic from 10:00 to 16:00
		\item a second peak of congestion in the afternoon, again triggered by congestion at Rondeau at around 16:00, spilling back on the whole freeway in around 60 minutes, and lasting approximatively two hours.
	\end{itemize}

\begin{figure}
\centering 
\includegraphics[width=0.4\textwidth]{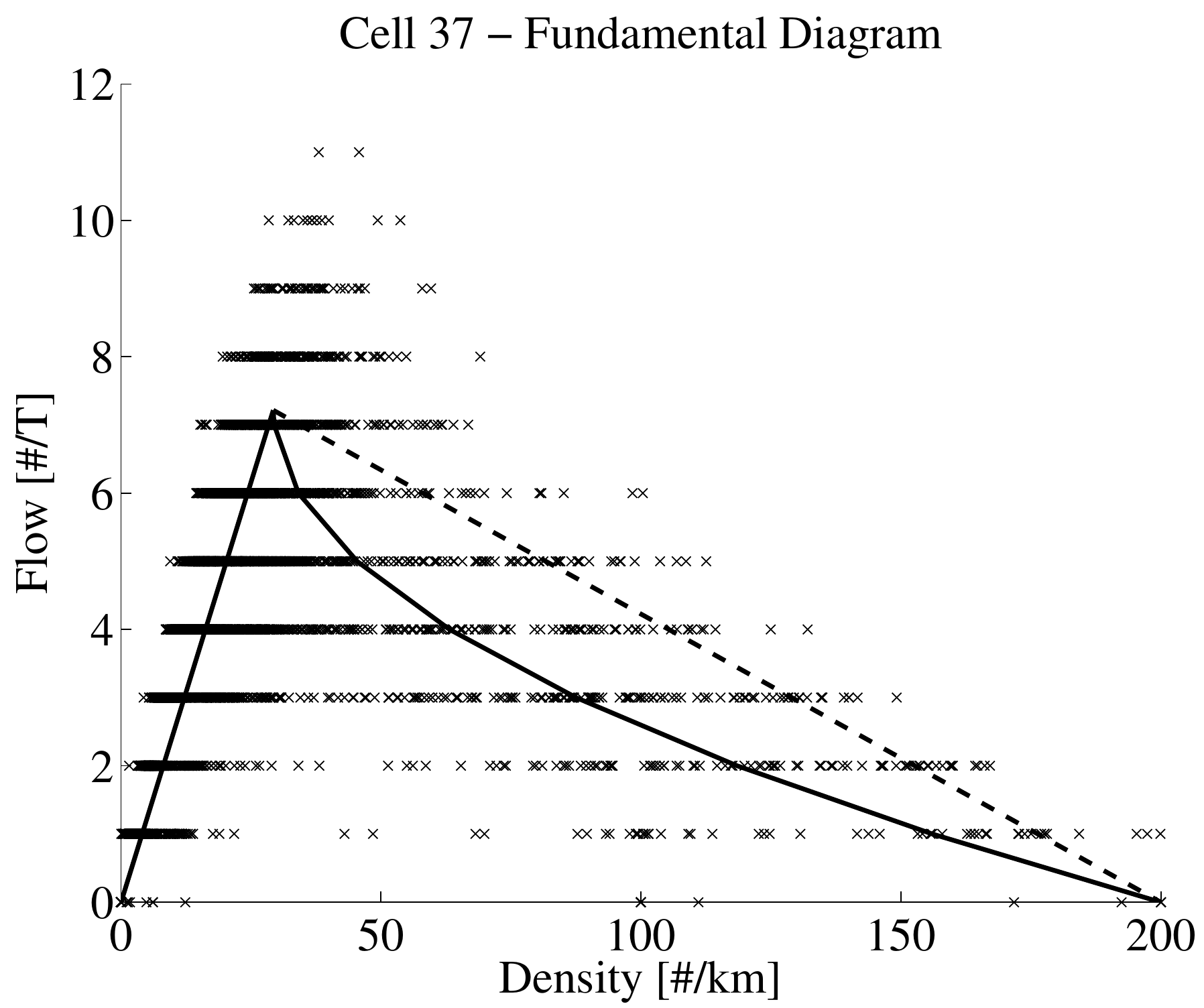}
\caption{Calibration of the Fundamental Diagram on the cell Eybens exit - slow lane. The linear-convex Fundamental Diagram, calibrated using data from April 10th, 2014, is shown in thick line. The dashed thick line represents the corresponding linear Fundamental Diagram in congested regime. Each cross is a (flow, density) pair measured on April 24th, 2014, one for each time slot of $T = 15$ seconds during the whole day. Flow is the number of counted vehicles crossing the sensing location, density is the measured density of vehicles, during the time slot.}
\label{figure:cell37_Eybensexitslow}
\end{figure}

As for the previous section, we consider the following simple rules to set the matrix of splitting ratios
	\begin{itemize}
		\item let $e$ be a fast lane cell. Then $70$\% of vehicles continue on the fast lane cell and $30$\% turn into the slow lane cell; if in the following section there are three parallel cells, vehicles split uniformly;
		\item let $e$ be a slow lane cell. If among the following cells there is not an offramp, then $70$\% of vehicles continue on the slow lane cell and $30$\% turn into the fast lane cell; if in the following section there are three parallel cells, vehicles split uniformly. Otherwise, $20\%$ of the flow is directed towards the offramp, and the rest splits as previously specified;
		\item if $e$ is an onramp cell and $j$ is the following slow ramp cell, then $R_{ej} = 1$.
		\item if $e$ is a queue cell and $j$ is the following onramp cell, then $R_{ej} = 1$.	
	\end{itemize}
In words, vehicles split according to a $70$\%-$30$\% lane-change rule in the cells on the main line, and at each offramp approximatively $10$\% of vehicles exit from the freeway.

\subsubsection{Implementation}
	
To assess our method, we considered the whole month of April 2014 (except April 13th, a Sunday, for which FCD measurements were not provided). A typical result of calibration of the Fundamental Diagram is illustrated in Figure~\ref{figure:cell37_Eybensexitslow}, which shows in thick black the linear-convex Fundamental Diagram, in dashed thick black the corresponding standard linear Fundamental Digram in congestion regime, and as crosses the pairs (density, flow) measured on a day different from that used for calibration, in this case April 24th, 2014. As standard and well known, data in freeflow regime are in good accordance with the linear part, while data in congested regime are much more scattered and more difficult to fit. As it can be noticed, the standard bilinear Fundamental Diagram overestimates the flows in congested regime (the dashed think line is on average higher than the corresponding pairs (density, flow)), while the convex quadratic curve seems to better capture the average flow-density relation. Nonetheless, it is clear that the so found curve is only a very crude approximation of such a relation, which might be better captured using a stochastic description \cite{JabariTRB:2013}. Investigation of the latter possibility will be the focus of future research.

The proposed algorithm was implemented in Matlab on a non dedicated commercial laptop with 2.1 GHz i7-4600U CPU and 8 GB RAM. Optimization problems, required both in offline and online steps, were solved using standard Matlab functions as well as the modelling and optimization system \texttt{CVX} \cite{CVX, GrantSpringer:08}. The time required for calibration of the Fundamental Diagrams is between $30$ and $40$ seconds for each cell, while reconstruction of all the samples for a whole day requires less $10$ minutes, averaging $100$ ms per sample. 

\begin{figure}
\begin{center}
\begin{tabular}{cc}
\includegraphics[width=0.45\textwidth]{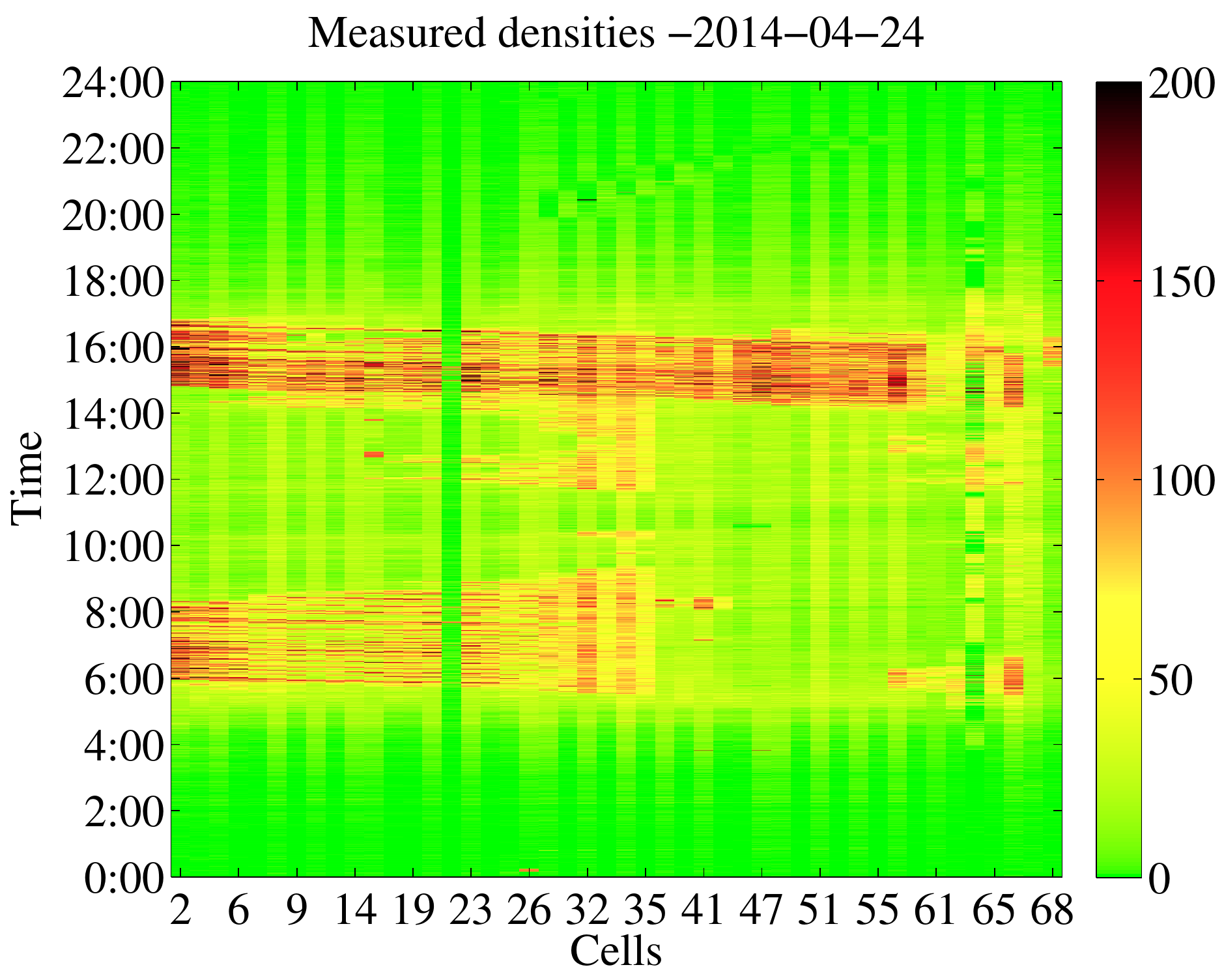}&
\includegraphics[width=0.45\textwidth]{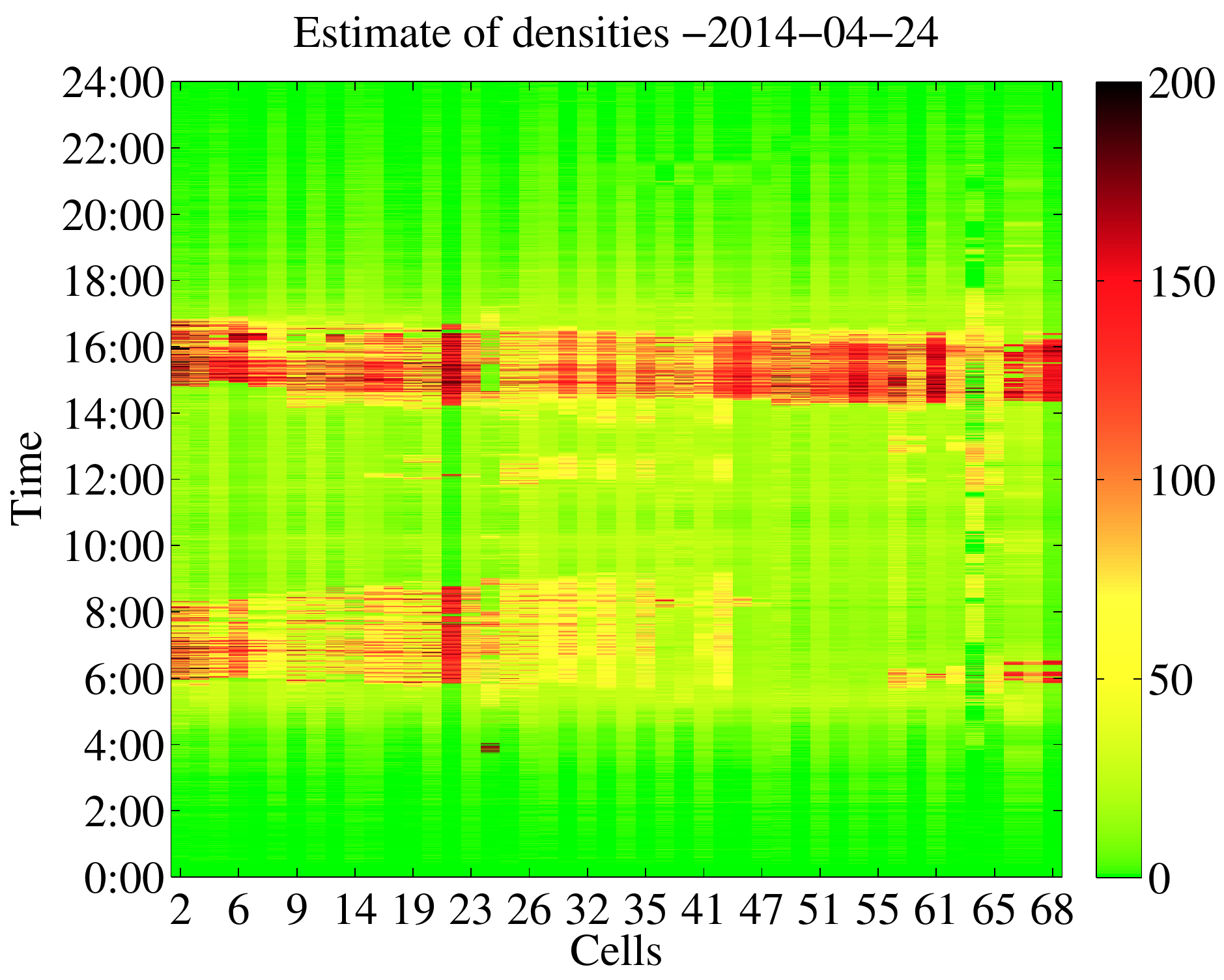}
\end{tabular}
\end{center}
\caption{Numerical results of density estimation on all cells: measured densities (left panel) and estimated densities right panel). Day: 24-04-2014}
\label{figure:NumericalExperiment_21042024_dens}
\end{figure}

\begin{figure}
\begin{center}
\begin{tabular}{cc}
\includegraphics[width=0.45\textwidth]{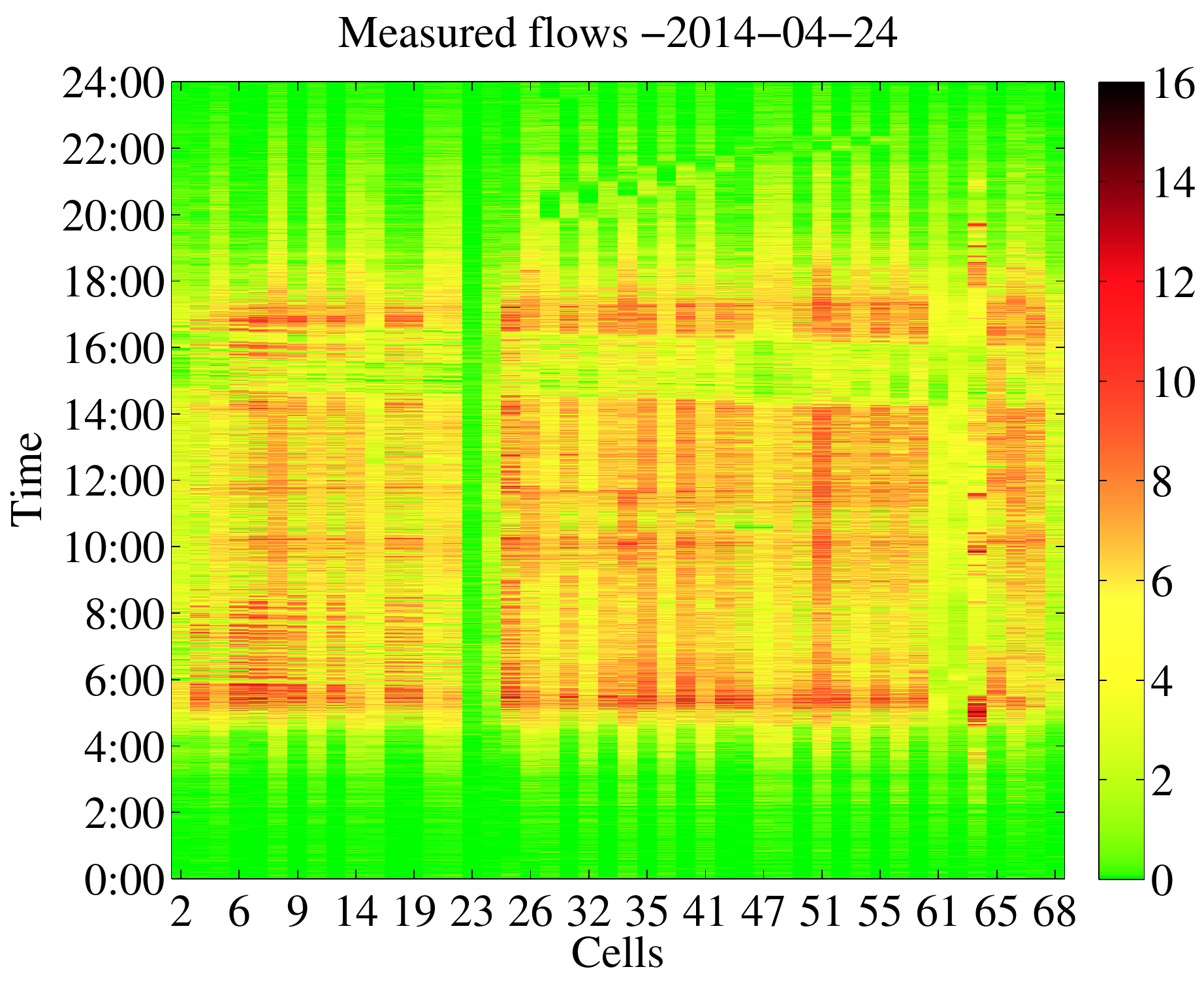}&
\includegraphics[width=0.45\textwidth]{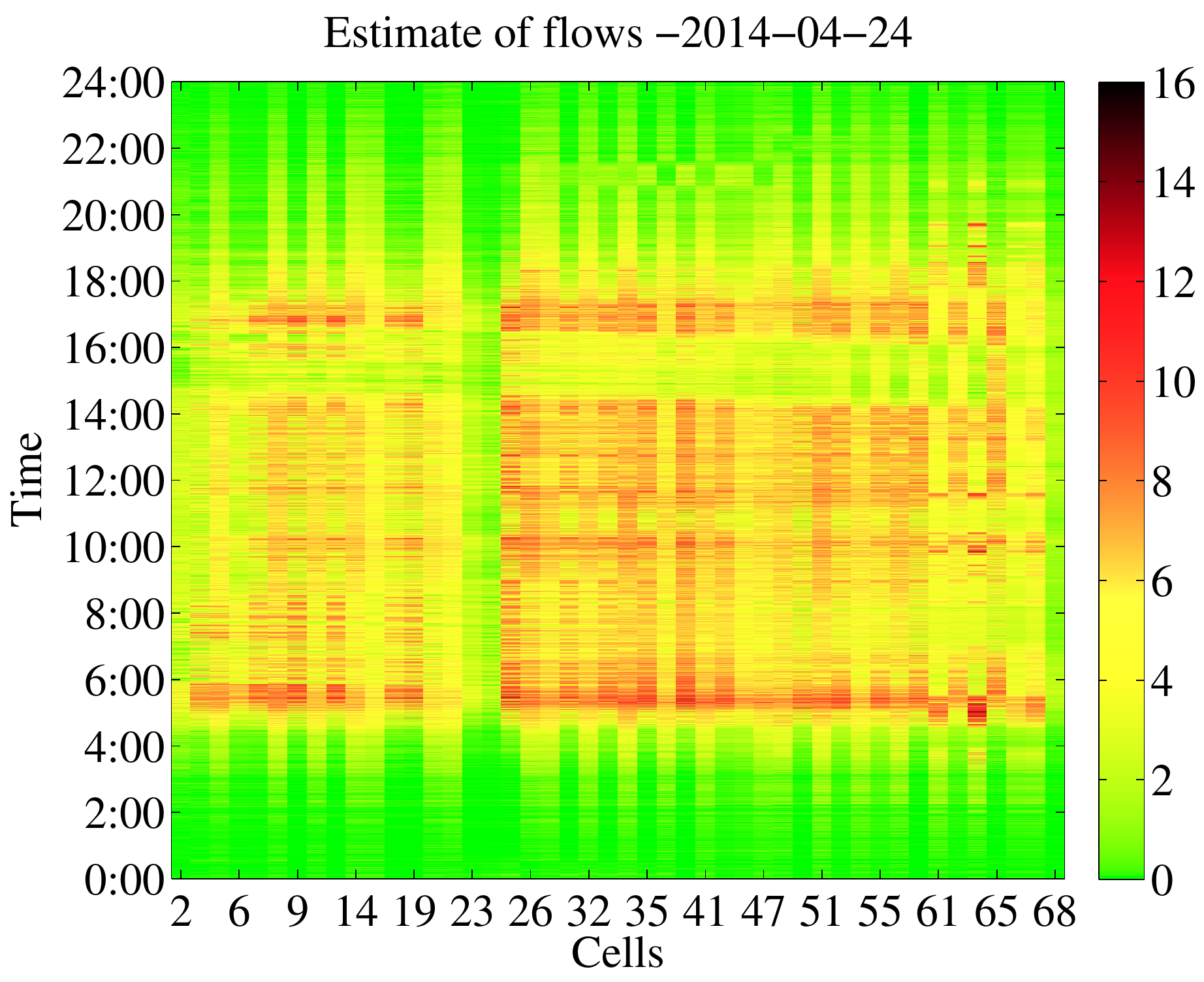}
\end{tabular}
\end{center}
\caption{Numerical results of flow estimation on all cells: measured densities (left panel) and estimated densities right panel). Day: 24-04-2014}
\label{figure:NumericalExperiment_21042024_flow}
\end{figure}

Typical results are reported in Figures~\ref{figure:NumericalExperiment_21042024_dens}-\ref{figure:NumericalExperiment_21042024_flow}. For validation purposes only, density and flow measurements from all GTL fixed sensors are considered ground truth. As such, the left panels show the evolution of the \emph{``true'' measured density and flow} in all the cells on the main line, over the whole day, for each day. On the $x$-axis, the $46$ sensing locations along the main line (numbers correspond to the labels in Figure~\ref{figure:RocadeSudStylized}), on the $y$-axis, the $5760$ time slots over the whole validation day (one slot every $T = 15$ seconds). The chosen colors range from green (low density/flow) to yellow (medium/critical density/flow) to red (high density/flow). In the right panels, we show, using the same legend, the results of density and flow reconstruction. As it can be observed, the estimation algorithm is able to represent the four congestion events described in the previous section in a reasonably good way, given the limited amount of information employed; in particular, observe that the two small congestions at Rondeau during the morning and at Eybens during early afternoon, when present, are both detected. Notice that the resulting estimate remains good and reasonably close to the real profile, despite the absence of flow measurements on ramps, which are not negligible, especially in peak hours.

\begin{table}
\begin{center}
\begin{tabular}{l||ccc||ccc}
						&	\multicolumn{3}{c}{$\delta:\,a^{\rho}(t,e)\leq \delta$} & \multicolumn{3}{c}{$\delta:\,a^{\varphi}(t,e)\leq \delta$} \\
Considered \%$(t,e)$:\	&	$75$\%		&	$90$\%		& $95$\%					&	$75$\%		&$90$\%						& $95$\%		\\
\hline
2014-04-01 -Tuesday &10.4317 &19.4079 &32.1121 &2.4901 &3.9824 &4.8838\\2014-04-02 -Wednesday &10.058 &16.8625 &23.4218 &2.504 &4.0508 &4.9457\\2014-04-03 -Thursday &14.214 &39.0677 &68.6383 &2.401 &3.8752 &4.765\\2014-04-04 -Friday &13.3274 &35.6493 &67.2034 &2.5477 &4.0648 &4.9543\\2014-04-05 -Saturday &7.6795 &11.6078 &14.9678 &2.5927 &4.138 &5.0691\\2014-04-06 -Sunday &5.6087 &7.8254 &9.5749 &2.1949 &3.5373 &4.4432\\2014-04-07 -Monday &8.8653 &14.9656 &21.5993 &2.2046 &3.6666 &4.6213\\2014-04-08 -Tuesday &9.4236 &16.1212 &24.1236 &2.134 &3.5462 &4.4686\\2014-04-09 -Wednesday &10.7057 &21.1416 &34.5579 &2.1496 &3.5922 &4.5308\\2014-04-10 -Thursday &10.3578 &21.7147 &38.2241 &2.1529 &3.5922 &4.5125\\2014-04-11 -Friday &10.1152 &19.0844 &32.2754 &1.8778 &3.2413 &4.1324\\2014-04-12 -Saturday &5.7334 &8.4439 &11.1728 &2.1883 &3.5934 &4.5191\\2014-04-14 -Monday &9.5789 &18.4657 &31.2646 &2.1349 &3.5506 &4.5182\\2014-04-15 -Tuesday &11.2095 &18.9568 &32.7132 &2.0675 &3.6966 &4.6661\\2014-04-16 -Wednesday &12.9045 &28.6587 &46.8904 &2.2881 &3.9458 &4.9609\\2014-04-17 -Thursday &17.482 &43.2213 &62.2287 &2.0069 &3.5265 &4.4383\\2014-04-18 -Friday &14.8889 &33.0979 &48.2054 &2.3525 &4.0964 &5.1424\\2014-04-19 -Saturday &6.3368 &9.8998 &12.8105 &2.3955 &3.9664 &4.8847\\2014-04-20 -Sunday &5.1875 &7.1362 &8.7217 &2.0231 &3.2701 &4.0776\\2014-04-21 -Monday &5.4528 &7.433 &8.9298 &2.1513 &3.4238 &4.2284\\2014-04-22 -Tuesday &12.0275 &26.3377 &43.2225 &2.2019 &3.7727 &4.7825\\2014-04-23 -Wednesday &12.9473 &32.5637 &49.0958 &2.0893 &3.4793 &4.4004\\2014-04-24 -Thursday &18.222 &42.5245 &59.2349 &2.0391 &3.3786 &4.282\\2014-04-25 -Friday &24.4905 &49.1946 &65.6698 &2.0058 &3.3379 &4.244\\2014-04-26 -Saturday &6.0991 &8.7496 &10.8339 &1.922 &3.3109 &4.1612\\2014-04-27 -Sunday &5.4655 &7.4201 &8.93 &2.1341 &3.3638 &4.1301\\2014-04-28 -Monday &13.8312 &35.9539 &57.5468 &1.9 &3.3551 &4.2986\\2014-04-29 -Tuesday &9.7221 &17.4268 &25.9634 &2.2904 &3.7966 &4.714\\2014-04-30 -Wednesday &12.736 &29.0402 &47.1253 &2.2473 &3.7883 &4.7363\\\hline Average&10.8656&22.3439&34.3882&2.1961&3.6531&4.5694

\end{tabular}
\end{center}
\caption{Quantitative measurement of the performance of the proposed algorithm. Maximum magnitude of the absolute error on densities (left) and flows (right) on $75$\%, $90$\% and $95$\% of the pairs (cell, time) between 07:00 and 19:00.}
\label{table:quantitativeResults}
\end{table}

The performance of the algorithm is quantitatively illustrated via the absolute error between measured and estimated flows and densities \eqref{eq:metrics_reconstruction}.

The results are reported in  Table~\ref{table:quantitativeResults}, in which we report the maximum absolute error $\delta$ on the $75$\%, $90$\% and $95$\% of the pairs (cell, time), for all cells on the main line and all samples during a day, for all the considered days, and for both densities and flows. The table shows an average error of less than $11$ veh/km for the $75$\% of pairs and less than $23$ veh/jm for the $90$\% of the pairs. In $5$\% of the (cell, time) pairs the error is higher but still less than around $35$ veh/km. As a comparison, we considered an oracle that knows
	\begin{itemize}
		\item the exact outflow $\fout_e(t)$ for \emph{every} cell $e$ and \emph{every} sample time $t$ (namely, for each $15$ seconds time slot);
		\item whether the cell is in freeflow or in congestion, for \emph{every} cell $e$ and \emph{every} sample time $t$.
	\end{itemize}
The average error obtained by the oracle is $5.8$ for the $75$\% of the pairs (cell, time) and $12.5$ for the $90$\% of the pairs, which are high even with the high amount of additional (and precise) information available to the oracle. Indeed, this confirms that estimation in traffic systems is a rather difficult task, and that errors of absolute magnitude around $10$-$20$ veh/km can be acceptable, as they capture the qualitative trend features of the traffic system - such as low density, medium density, high density. On the other side, estimate of flows are rather low, being less than $1\sim 2$ vehicles for the $75$\% of the pairs (cell, time), and less than $2\sim 3$ for the $90$\% of the pairs.

The biggest difference between estimated and measured densities can be observed at Eybens exit (cells 37-38), where the estimated flow constantly predicts a higher density than the measured one. The explanation is however very straightforward: as mentioned, we do not use ramp data in order to show the prowess of our method even employing a very small number of sensors. On the other hand, as mentioned in the description of the data, the exit of Eybens is a critical point whose ramp is selected by a high fraction of vehicles to exit the Rocade towards the town. However, the corresponding cells belong to a long FCD segment running from Saint-Martin-d'H\`{e}res (cells 32-33) to Eybens entrance (cells 41-42), which provides, during peak time, just one set of rather low speed measurements, which do not distinguish between the stretch of road before Eybens exit (congested and at low speed) against that after Eybens exit (uncongested and at high speed). Due to the so obtained low speed measurement, the algorithm tends to estimate a high number of vehicles along the whole segment, instead of two different regimes before and after the ramp. Analogously, we observe a mismatch between measured and estimated flow in the first and last sections, which are due to unobserved flows from onramp and to offramps. A second discrepancy is the smoothness of the reconstructed density and flow, as compared with the more scattered measurements. The latter is due to the high measurement rate, which during stop-and-go phenomena results in measurements which rapidly oscillate between stopped vehicles and low or medium speed. On the converse, the optimization based flow reconstruction and the first order mass conservation law for densities have a low pass effect therefore producing smoother, more regular patterns. Further research direction will investigate the possibility to detect stop-and-go phenomena and reproduce, at least qualitatively, the resulting irregular patterns.

\section{Conclusions}
\label{sec:conclusions}

This paper addressed the problems of data fusion of heterogeneous sources of information for density estimation in Road Transportation Networks and optimal sensor placement via a heuristic that we called Virtual Variance algorithm. A gradient descent procedure for the calibration of the Fundamental Diagram is also discussed. Efficacy of the proposed solutions is illustrated on a regular grid and on the real world scenario of the Rocade Sud in Grenoble. Future research directions include and are not limited to estimation of statistical properties of measurement noises from real data and development of stochastic models for the relation between flows, speed and densities, optimization of the observer's parameters for minimization of mean-square reconstruction error, calibration of the matrix of splitting ratios, and extension of the optimal sensor placement strategy to maximize density reconstruction performance.

\section{Acknowledgements}
\label{sec:ack}

The authors gratefully thank Rene Fritz, Mike Corlett, and INRIX Europe for providing the Floating Car Data that led to the results presented in this work. We also thank Roland Dollet and Sylvain Nachef for introducing us to the topic of Density Reconstruction in Road Transportation Systems.

\bibliographystyle{unsrt}
\bibliography{bibliography}

\end{document}